\documentclass{amsproc}
\usepackage{amsmath,amsfonts,amssymb,amsthm}
\title{Examples of noncommutative instantons}

\author{Giovanni Landi}
\address{Dipartimento di Matematica e Informatica, Universit\`a di Trieste\\
Via A. Valerio 12/1, I-34127 Trieste, Italy}
\email{landi@univ.trieste.it}


\newtheorem{theo}{Theorem}[section]
\newtheorem{lma}[theo]{Lemma}
\newtheorem{prop}[theo]{Proposition}
\newtheorem{coro}[theo]{Corollary}

\theoremstyle{definition}
\newtheorem{defi}[theo]{Definition}
\newtheorem{exam}[theo]{Example}

\theoremstyle{remark}
\newtheorem{rema}[theo]{Remark}

\numberwithin{equation}{section}

\renewcommand{\bar}[1]{\overline{#1}}
\renewcommand{\tilde}[1]{\widetilde{#1}}

\newcommand{\half}{{\mathchoice{\oh}{\oh}{\shalf}{\shalf}}} 
\newcommand{\oh}{{\tfrac{1}{2}}}    
\newcommand{\shalf}{{\scriptstyle\frac{1}{2}}} 
\def\II{\mathbb{I}}
\def\ii{\mathrm{i}} 
\def\into{\hookrightarrow}

\def\isom{\simeq}

\def\C{\mathbb{C}} 

\def\I{\mathbb{I}} 
 
\def\bR{\mathbb{R}}

\def\bT{\mathbb{T}}
 
\def\A{\mathcal{A}}

\def\B{\mathcal{B}}

\def\E{\mathcal{E}}

\def\cH{\mathcal{H}}

\def\cO{\mathcal{O}}

\def\cU{\mathcal{U}}

\newcommand{\ket}[1]{|#1\rangle}    
\newcommand{\bra}[1]{\langle#1|}    

\def\chern{\mathrm{ch}}

\def\dix{\mathrm{Tr}_\omega}
\def\ind{\mathrm{index}\mb}

\def\ad{\mathrm{ad}}

\def\Cinf{C^\infty}
\def\class{{(0)}}

\def\M{M_\theta}

\def\S{S_\theta}
\def\Sk{S_{\theta'}}

\def\Top{\mathrm{Top}}

\def\YM{\mathrm{YM}}
 
\def\lt{\triangleright}    
\def\rt{\triangleleft}    

\def\be{\begin{equation}}
\def\ee{\end{equation}}
\def\bea{\begin{eqnarray}}
\def\eea{\end{eqnarray}}
\def\bmult{\begin{multline}}
\def\emult{\end{multline}}
\def\nn{\nonumber}
\def\mb{\mbox{ }}

\newcommand{\bean}{\begin{eqnarray*}}
\newcommand{\eean}{\end{eqnarray*}}
\newcommand{\cinf}{C^\infty}       
\newcommand{\ca}{{\mathcal A}}

\newcommand{\cc}{{\mathcal C}}

\newcommand{\ce}{{\mathcal E}}
\newcommand{\cf}{{\mathcal F}}

\newcommand{\ch}{{\mathcal H}}

\newcommand{\ck}{{\mathcal K}}

\newcommand{\cs}{{\mathcal S}}

\newcommand{\cu}{{\mathcal U}}

\newcommand{\IC}{{\mathbb C}}

\newcommand{\IN}{{\mathbb N}}

\newcommand{\IR}{{\mathbb R}}

\newcommand{\IT}{{\mathbb T}}
\newcommand{\IZ}{{\mathbb Z}}

\def\lra{\longrightarrow}

\newcommand{\wh}{\widehat}

\def\bar#1{\overline{#1}}

\newcommand{\omca}{\Omega {\mathcal A}}

\newcommand{\oca}[1]{\Omega^{#1}{\mathcal A}}

\newcommand{\comca}{{\mathcal E}\otimes_{\mathcal A}\Omega {\mathcal A}}
\newcommand{\coca}[1]{{\mathcal E}\otimes_{\mathcal A}\Omega^{#1}{\mathcal A}}
\newcommand{\cha}[1]{C_{#1}({\mathcal A})}
\newcommand{\ccha}[1]{CC_{#1}({\mathcal A})}

\newcommand{\ota}{\otimes_{\mathcal A}}

\newcommand{\otc}{\otimes_{\IC}}
\newcommand{\ot}{\otimes}

\newcommand{\op}{\oplus}

\newcommand{\pot}{\overset{.}{\otimes}}

\newcommand{\lr}{\longrightarrow}
\newcommand{\uno}{x_{\mbox{\tiny{1}}}}
\newcommand{\unob}{\bar{x}^{\mbox{\tiny{1}}}}
\newcommand{\due}{x_{\mbox{\tiny{2}}}}
\newcommand{\dueb}{\bar{x}^{\mbox{\tiny{2}}}}
\newcommand{\tre}{x_{\mbox{\tiny{3}}}}
\newcommand{\treb}{\bar{x}^{\mbox{\tiny{3}}}}
\newcommand{\qu}{x_{\mbox{\tiny{4}}}}
\newcommand{\qub}{\bar{x}^{\mbox{\tiny{4}}}}

\newcommand{\hil}{\ensuremath{\mathcal{H}}}

\newcommand{\hs}[2]{\left\langle #1,#2\right\rangle}

\DeclareMathOperator{\id}{id}
\DeclareMathOperator{\Mat}{Mat}
\DeclareMathOperator{\Tr}{Tr}

\DeclareMathOperator{\tr}{tr}
\DeclareMathOperator{\SU}{SU}

\DeclareMathOperator{\U}{U} 
\DeclareMathOperator{\SO}{SO}
\DeclareMathOperator{\Spin}{Spin}
 
\DeclareMathOperator{\Hom}{Hom}
\DeclareMathOperator{\End}{End}
\newcommand{\dd}{{\rm d}} 
\newcommand{\btheo}{\begin{theo} ~~\\}
\newcommand{\etheo}{\hfill $\Box$  \end{theo} }
\newcommand{\bprop}{\begin{prop} ~~\\ }
\newcommand{\eprop}{\hfill $\Box$ \end{prop} }
\newcommand{\bdefi}{\begin{defi} ~~\\ \rm }
\newcommand{\edefi}{\hfill $\Box$ \end{defi} }
\newcommand{\brema}{\begin{rema} ~~\\ }
\newcommand{\erema}{ \hfill $\Box$ \end{rema} }
\newcommand{\bexam}{\begin{exam} ~~\\ }
\newcommand{\eexam}{\hfill $\Box$ \end{exam} }
\newbox\ncintdbox \newbox\ncinttbox
\setbox0=\hbox{$-$} \setbox2=\hbox{$\displaystyle\int$}
\setbox\ncintdbox=\hbox{\rlap{\hbox to
\wd2{\hskip-.125em\box2\relax \hfil}}\box0\kern.1em}
\setbox0=\hbox{$\vcenter{\hrule width 4pt}$}
\setbox2=\hbox{$\textstyle\int$}
\setbox\ncinttbox=\hbox{\rlap{\hbox to
\wd2{\hskip-.175em\box2\relax \hfil}}\box0\kern.1em}
\newcommand{\ncint}{\mathop{\mathchoice{\copy\ncintdbox}%
{\copy\ncinttbox}{\copy\ncinttbox}{\copy\ncinttbox}}\nolimits}

\begin{document}

\maketitle


\tableofcontents

\date{November 24, 2006}

\section{Introduction}

These notes aim at a pedagogical introduction to recent work on deformation of spaces and deformation of vector bundles over them, which are relevant both in mathematics and in physics, notably monopole and instanton  bundles. 

We first decribe toric noncommutative manifolds (also known as isospectral deformations).
These come from deforming the usual Riemannian geometry of a manifold along a torus embedded in the isometry group, thus obtaining a family of noncommutative geometries which are isospectral to the starting one.

To be specific, we give a detailed introduction to gauge theories on a toric four-sphere $\S^4$, including the definition of a Yang--Mills action functional with associated equations of motion and self-duality equations. Solutions of the latter equations -- that is  {\it instantons} -- are absolute minima of the action functional. A particular class of them is construct from a $\SU(2)$ principal bundle $\Sk^7 \rightarrow \S^4$, with a suitable use of twisted conformal symmetries yielding a five parameter family of instantons.  

In the second part, we describe a different quantum principal bundle $S_q^7 \rightarrow S_q^4$ having as  `structure group' the quantum group $\SU_q(2)$. The quantum sphere $S^7_q$  arises from the symplectic quantum group $Sp_q(2)$ and a quantum $4$-sphere $S^4_q$ is obtained via a suitable self-adjoint idempotent $p$ whose entries generate the algebra $A(S^4_q)$ of polynomial functions over it -- a property in common with the toric four sphere $\S^4$. This projection determines a deformation of an instanton bundle over the classical sphere $S^4$.
For this bundle a construction of a Yang--Mills action functional  and self-duality equations is still missing. 

To be definite, we think of an instanton as  `a rank two complex vector bundle on a four dimensional manifold endowed with a self-dual connection'  but, of course, we do not mind proper generalizations. 
In  this respect, we mention the work \cite{NS98} on noncommutative instantons on a noncommutative $\IR^4$ which has resulted in a host of interesting and current developing activities.

\subsection*{Acknowledgment}
I wish to thank the organizers of the 2005 Summer School on ``Geometric and Topological  Methods for Quantum  Field Theory'', July 11-29 2005, Villa de Leyva, Colombia, for the kind invitation to deliver these lectures.
It is a great pleasure to remember Sergio Adarve, Alexander Cardona, Hern\'an Ocampo, Sylvie Paycha, Bernardo Uribe and all other peoples at Villa de Leyva -- with a special mention for Marta Kovacsics --
for the fantastic time in the Colombian Andes. These notes are based on the papers \cite{LPR06,LS04,LS06}.

\section{Toric noncommutative manifolds}\label{se:toric-ncm}

Toric noncommutative manifolds are deformations of a classical Riemannian manifold and  satisfy all the properties  of a noncommutative spin geometry~\cite{C96}. They are the content of the following result taken from~\cite{CL01}, where they were called isospectral deformations. 

\begin{theo}\label{Theorem6}
Let $M$ be a compact spin Riemannian manifold whose isometry group 
has rank $r \geq 2$. Then $M$ admits a natural 
one parameter isospectral deformation to noncommutative geometries $M_{\theta}$.
\end{theo}

The idea is to deform the standard spectral triple describing the Riemannian geometry of $M$ along a torus embedded in the isometry group of $M$, thus obtaining a family of spectral
triples describing noncommutative geometries which are isospectral to the starting one of $M$.

\subsection{Deforming along a torus}\label{subse:def-torus}
Let $M$ be an $m$ dimensional compact Riemannian manifold equipped with an isometric smooth action $\sigma$ of an $n$-torus $\bT^n$, with $n \geq 2$.  We denote by $\sigma$ also the corresponding action of $\bT^n$ by automorphisms on the algebra $\Cinf(M)$, obtained by pull-back. The algebra $\cinf(M)$ may be decomposed
into spectral subspaces which are indexed by the dual group $\IZ^n = \wh\IT^n$. Now, with $s=(s_1, \cdots s_n)\in\bT^n$, 
each $r \in \IZ^n$ labels a character $e^{2\pi i s} \mapsto e^{2\pi i r\cdot s}$ of $\IT^n$, with scalar product
$r\cdot s := \sum_{j=1}^{n} r_j s_j$. The $r$-th spectral
subspace for the action $\sigma$ of $\IT^n$ on $\cinf(M)$ consists of those smooth
functions $f_r$ for which
\be\label{actor}
\sigma_s (f_r) = e^{2\pi \ii r\cdot s} \,f_r ,
\ee
and each $f \in \cinf(M)$ is the sum of a unique rapidly convergent series $f = \sum_{r\in\IZ^n} f_r$.
Let now $\theta = (\theta_{j k} = - \theta_{k j})$ be a real antisymmetric $n\times n$ matrix. 
The $\theta$-deformation of $\cinf(M)$ may be defined by replacing the ordinary product by a
deformed product given on spectral subspaces by
\be
\label{eq:star-product}
f_r \times_\theta g_{r'} := f_r ~ \sigma_{\half r \cdot \theta}( g_{r'}) = 
e^{  \pi \ii r \cdot \theta \cdot r' } f_r g_{r'} ,
\ee
where $r\cdot \theta =(r_j \theta_{j 1}, \ldots, r_j \theta_{j n} )\in \bT^n$. This product is then extended linearly to all functions in $\Cinf(M)$. We denote $\cinf(M_\theta) := (\cinf(M),\times_\theta)$ and note that the action $\sigma$ of $\IT^n$ on $\cinf(M)$ extends to an action on $\cinf(M_\theta)$ given again by \eqref{actor} on the homogeneous elements.

Next, let us take $M$ to be a spin manifold with $\cH:=L^2(M,\cs)$ the Hilbert space of spinors and $D$ the usual Dirac operator of the metric of $M$. Smooth functions act on spinors by pointwise multiplication thus giving a representation $\pi : \Cinf(M) \to \B(\cH)$, the latter being the algebra of bounded operators on $\ch$. 

There is a double cover $c: \widetilde{\bT}^n \to \bT^n$ and a representation of $\tilde \bT^n$ on $\cH$ by unitary operators $U(s), s \in \tilde\bT^n$, so that 
\[
U(s) D U(s)^{-1} = D,
\]
since the torus action is assumed to be isometric, and such that for all $f \in \Cinf(M)$, 
\[
U(s) \pi(f) U(s)^{-1} = \pi(\sigma_{c(s)}(f)).
\]
Recall that an element $T\in \B(\cH)$ is called smooth for the action of $\tilde \bT^n$ if the map 
\[
\tilde \bT^n \ni s \mapsto \alpha_s(T) := U(s) T U(s)^{-1},
\]
is smooth for the norm topology. From its very definition, the map $\alpha_s$ coincides on $\pi(C^\infty(M)) \subset \B(\cH)$ with the automorphism $\sigma_{c(s)}$. Much as it was done before for the smooth functions, we shall use the torus action to give a spectral decomposition of smooth elements of $\B(\cH)$.  Any such a smooth element $T$ is written as a (rapidly convergent) series $T =\sum T_{r}$ with $r\in\IZ$ and each $T_{r}$ is homogeneous of degree $r$ under the action of $\tilde \bT^n$, {i.e.}
\be\label{homocomp}
\alpha_s(T_{r}) =e^{2 \pi \ii r \cdot s } T_{r} ,  \quad \forall \quad s  \in \tilde \bT^n .
\ee
Let $p=(p_1, p_2,\ldots, p_n)$ be the infinitesimal generators of the action of $\tilde \bT^n$ so that we can write $U(s)=\exp{2 \pi \ii s \cdot p}$.  
Now, with $\theta$ a real $n\times n$ anti-symmetric matrix as above,  a twisted representation 
of the smooth elements $\B(\cH)$ on $\cH$ is defined by
\be\label{twist}
L_\theta(T):=\sum_r T_r U\big(\half( r_j \theta_{j1}, \ldots,  r_j\theta_{jn})\big) 
= \sum_r T_r \exp \big( \pi \ii \,  r_j \theta_{jk} p_k \big), 
\ee
Taking smooth functions on $M$ as elements of $\B(\cH)$, via the representation $\pi$, the previous 
definition gives an algebra $L_\theta(\Cinf(M))$ which we may think of as a representation (as bounded operators on $\ch$) of the algebra $\Cinf(\M)$. Indeed, by the very definition of the product $\times_\theta$ in \eqref{eq:star-product} one establishes  that 
\[
L_\theta(f\times_\theta g)= L_\theta(f) L_\theta(g),
\]
proving that the algebra $\Cinf(M)$ equipped with the product $\times_\theta$ is isomorphic to the algebra $L_\theta(\Cinf(M))$. We shall think of $L_\theta$ as a {\it quantization map}
\[
L_\theta: \Cinf(M) \to \Cinf(\M).
\]
This quantization map will play a key role in what follows, and allows one to extend differential geometric properties from $M$ to the noncommutative space $\M$.

It is shown in \cite{Ri93} that there is a natural completion of the algebra $\cinf(M_\theta)$ to a $C^*$-algebra $C(M_\theta)$ whose smooth subalgebra -- under the extended action of $\IT^n$ -- is
precisely $\cinf(M_\theta)$. Thus, we can understand $L_\theta$ as a quantization map from $\Cinf(M)$ to $\Cinf(\M)$ providing a strict deformation quantization in the sense of Rieffel.
More generally, in \cite{Ri93} one considers a (not necessarily commutative) $C^*$-algebra $A$ carrying an action of $\bR^n$. For an anti-symmetric $n \times n$ matrix $\theta$, one defines a deformed product $\times_\theta$ between elements in $A$ much as we did before. The algebra $A$ equipped with the product $\times_\theta$ gives rise to a $C^*$-algebra denoted by $A_\theta$. Then the collection $\{A_{\hbar\theta} \}_{\hbar\in [0,1]}$ is a continuous family of $C^*$-algebras providing a strict deformation quantization in the direction of the Poisson structure on $A$ defined by the matrix $\theta$. 
Our case of interest corresponds to the choice $A=C(M)$ with an action of $\bR^n$ that is periodic or, in other words, an action of $\bT^n$. The smooth elements in  the deformed algebra make up
the algebra $\Cinf(\M)$.

It was shown in \cite{CL01} that the triple $(L_\theta(C^\infty(M)), \cH, D)$ satisfies  all axioms of a noncommutative spin geometry  \cite{C96}. There is also a grading $\gamma$ (for the even case) and a real structure $J$. 
In particular, boundedness of the commutators $[D,L_\theta(f)]$ for $f \in C^\infty(M)$ follows from 
$[D,L_\theta(f)]=L_\theta([D,f])$, $D$ being of degree $0$ (since $\bT^n$ acts by isometries, each infinitesimal generator $p_k$ commutes with $D$). This noncommutative geometry is an isospectral deformation of the classical Riemannian geometry of $M$, in that the spectrum of the operator $D$ coincides with that of the Dirac operator $D$ on $M$. 
Thus all spectral properties are unchanged. In particular, the triples are $m^+$-summable and there is a noncommutative integral as a Dixmier trace \cite{Dix66},
\be\label{dix}
\ncint L_\theta (f) := \dix \big( L_\theta (f) |D|^{-m}  \big), 
\ee
with $f \in \Cinf(\M)$ understood in its representation on $\cH$. A drastic simplification~\cite{GIV05} of this noncommutative integral is given by the  
\begin{lma} \label{lma:dix}
If $f \in \Cinf(M)$ then
$$\ncint L_\theta (f) = \int_{M} f \dd \nu.$$
\end{lma}

\bigskip

Vector bundles on $\M$ were described in \cite{LS06} to which we refer for details and proofs. 
Crucially, they can be given in terms of a deformed product and action.

Let $E$ be a {\it $\sigma$-equivariant} vector bundle  $M$, that is a bundle which carries an action $V$ of $\bT^n$ by automorphisms, covering the action {\it $\sigma$} of $\bT^n$ on $M$,
\be
\label{eq:sigma-equivariant}
V_s (f \psi) = \sigma_s(f) V_s(\psi), \quad \forall \, f \in \Cinf(M) , \, \psi \in \Gamma(M,E) .
\ee
The $\Cinf(\M)$-bimodule $\Gamma(\M, E)$ is defined as the vector space $\Gamma(M,E)$,
that is we have a quantization map
\[
L_\theta : \Gamma(\M, E) \to \Gamma(M,E), 
\]
but with a twisted bimodule structure given by 
\be\label{eq:bi-action}
f \lt_\theta \psi = \sum_k f_k V_{\half k\cdot \theta} (\psi) , \qquad
\psi \rt_\theta f = \sum_k V_{-\half k\cdot \theta} (\psi) f_k ,
\ee
where $f=\sum_k f_k$ with $f_k \in \Cinf(M)$ homogeneous of degree $k$ under the action of $\bT^n$ -- as in \eqref{homocomp} -- and $\psi$ is a smooth section of $E$. That these are indeed actions of $\Cinf(\M)$ can be estalished  with the use of the explicit expression \eqref{eq:star-product} for the deformed product  and of 
equation~\eqref{eq:sigma-equivariant} for the equivariance. 

The $\Cinf(\M)$-bimodule $\Gamma(\M,E)$ is finite projective \cite{CD02} and still carries an action of $\bT^n$ by $V$ with equivariance as in equation \eqref{eq:sigma-equivariant} for both the left and right action of $\Cinf(\M)$. Indeed, the group $\bT^n$ being abelian, one establishes that 
\[
V_s (f \lt_\theta \psi) = \sigma_s(f) \lt_\theta V_s(\psi), \quad \forall \, f \in \Cinf(\M) , \, \psi \in \Gamma(\M,E),
\]
and a similar property for the right structure $\rt_\theta$.

Although the left and right actions in \eqref{eq:bi-action} were defined with respect to an action of $\bT^n$ on $E$, the same construction can be done for vector bundles carrying an action of the double cover $\tilde \bT^n$. We have already seen an instance of this for the spinor bundle, where we defined a left action of $\Cinf(\M)$ on the spinors using \eqref{twist}. 

\subsection{Differential calculus on $\M$}\label{se:diff-calc}

A differential calculus on $\M$ can be constructed in two equivalent manners, either by extending to forms the quantization maps, or by using the general construction via the Dirac operator \cite{C94}. 

Firstly, let $(\Omega(M),\dd)$ be the usual differential calculus on $M$, with $\dd$ the exterior derivative. The quantization map $L_\theta: \Cinf(M) \to \Cinf(\M)$ is extended to $\Omega(M)$ by imposing that it commutes with $\dd$. The image $L_\theta(\Omega(M))$ will be denoted $\Omega(\M)$. 
Equivalently, $\Omega(\M)$ could be defined to be $\Omega(M)$ as a vector space but equipped with an `exterior deformed product' which is the extension of the product  \eqref{eq:star-product}  to $\Omega(M)$ by the requirement  that it commutes with $\dd$. 
Indeed, since the action of $\bT^n$ commutes with $\dd$, an element in $\Omega(M)$ can be decomposed  into a sum of homogeneous elements for the action of $\bT^n$ -- as was done for $\Cinf(M)$. Then one defines a deformed product $\times_\theta$ on homogeneous elements in $\Omega(M)$ as in \eqref{eq:star-product} and denotes $\Omega(\M)=(\Omega(M), \times_\theta)$. This construction is in concordance with the previous section, when $\Omega(M)$ is considered as the $\Cinf(M)$-bimodule of sections of the cotangent bundle. The extended action  
of $\bT^n$  from $\Cinf(M)$ to $\Omega(M)$ is used to endow the space $\Omega(\M)$ with the structure of a $\Cinf(\M)$-bimodule with the left and right actions given as in \eqref{eq:bi-action}. 

Next, we describe the differential calculus $\Omega_D(\Cinf(\M))$ on $\Cinf(\M)$ obtained from the Dirac operator $D$ of the isospectral spin geometry on $\M$.
Elements of the $\Cinf(\M)$-bimodule $\Omega^p_D(\Cinf(\M))$ of $p$-forms are classes of operators of the type
\[
\omega = \sum_j a_0^j [D, a_1^j]\cdots [D, a_p^j], \quad a_i^j \in \Cinf(\M), 
\]
modulo the sub-bimodule of operators 
\[
\big\{ \sum_j [D,b_0^j] [D, b_1^j]\cdots [D, b_{p-1}^j] ;~ 
\sum_j b_0^j  [D, b_1^j]\cdots [D, b_{p-1}^j]=0 \big\},
\]
with $b_i^j \in \Cinf(\M)$. The exterior differential $\dd_D$ is given by 
\[
\dd_D\bigg[\sum_j a_0^j [D, a_1^j]\cdots [D, a_p^j]\bigg] = \bigg[\sum_j [D,a_0^j] [D, a_1^j]\cdots [D, a_p^j]\bigg],
\]
and satisfies $\dd_D^2=0$. An inner product on forms is defined by declaring that forms of different degree are orthogonal, while for two $p$-forms $\omega_1,\omega_2$, it is given by
\be\label{dix2}
(\omega_1,\omega_2 )_D=\ncint \omega_1^* \omega_2.
\ee
Here the noncommutative integral is the natural extension of the one in \eqref{dix},
\[
\ncint T  := \dix \big( T |D|^{-m}  \big), 
\]
with $T$ an element in a suitable (summability) class of operators.
Not surprisingly, these two construction of forms agree: the differential calculi $\Omega(\M)$ and $\Omega_D(\Cinf(\M))$ are isomorphic \cite{CD02}. In particular, this allows one to integrate forms of top dimension, by defining
\[
\int_{\M} \omega := \ncint \omega_D, \quad \omega \in \Omega(\M),
\]
where $\omega_D$ denotes the element in $\Omega_D(\Cinf(\M))$ corresponding to $\omega$ (replacing every $\dd$ in $\omega$ by $\dd_D$). There is a noncommutative Stokes theorem, 
\begin{lma}
\label{lma:stokes}
If $\omega \in \Omega^{\dim M-1} (\M)$ then
$$
\int_{\M} \dd \omega = 0 .
$$
\end{lma}
The next ingredient is a Hodge star operator on $\Omega(\M)$. Classically, the Hodge star operator is a map $\ast: \Omega^p(M) \to \Omega^{m-p}(M)$ depending only on the conformal class of the metric on $M$. On the one end, since $\bT^n$ acts by isometries, it leaves the conformal structure invariant and therefore, it commutes with $\ast$. On the other hand, with the isospectral deformation one does not change the metric.  
Thus it makes sense to define a map $\ast_\theta : \Omega^p(\M)\to \Omega^{m-p}(\M)$ by 
\be
\ast_\theta L_\theta(\omega) = L_\theta(\ast \omega), \quad \mathrm{for} \quad \omega \in \Omega(\M) .
\ee

With this Hodge operator, there is an alternative definition of an inner product on $\Omega(\M)$. 
Given that $\ast_\theta$ maps $\Omega^p(\M)$ to $\Omega^{m-p}(\M)$, we can define 
\be\label{eq:inner-product-forms}
( \alpha,\beta )_2 = \ncint \ast_\theta(\alpha^* \ast_\theta \beta),    \forall \, \alpha,\beta \in \Omega^p(\M),
\ee
since $\ast_\theta(\alpha^* \ast_\theta \beta)$ is an element in $\Cinf(\M)$. Needless to say, under the isomorphism $\Omega_D(\Cinf(\M)) \simeq \Omega(\M)$, the inner product $(\cdot,\cdot)_2$ coincides with $(\cdot,\cdot)_D$ given in \eqref{dix2}.
\begin{lma}
\label{lma:d-dstar}
The formal adjoint $\dd^*$ of $\dd$ with respect to the inner product $(\cdot,\cdot)_2$, {i.e.} so that $(\dd^* \alpha, \beta)_2=(\alpha,\dd \beta)_2$,  is given on $\Omega^p(\M)$ by 
$$
\dd^* = (-1)^{m(p+1)+1} \ast_\theta \dd \ast_\theta .
$$
\end{lma}

\section{Gauge theory on the sphere $\S^4$}\label{se:gts}

We decribe noncommutative gauge theory with the specific model of an $\SU(2)$ noncommutative principal bundle $\Sk^7 \to \S^4$. This will also work as a general scheme for any $\theta$-deformed $G$-principal bundle (with $G$ a compact semisimple Lie group). 

\subsection{The principal fibration $\Sk^7 \to \S^4$}
The $\SU(2)$ noncommutative principal fibration $\Sk^7 \to \S^4$ is given by an algebra inclusion
$\A(\S^4) \into \A(\Sk^7)$. 
With $\theta$ a real parameter and deformation parameters given by
\[
\lambda_{1 2} = \bar{\lambda}_{2 1} =: \lambda=e^{2\pi \ii \theta}, 
\quad \lambda_{j 0} = \lambda_{0 j } = 1, \quad j=1,2 ,
\]
the algebra $A(\S^4)$ of polynomial functions on the sphere $\S^4$ is 
 generated by elements  $z_0=z_0^*, z_j, z_j^*$, $j=1,2$, subject to relations
\be\label{s4t}
z_\mu z_\nu = \lambda_{\mu\nu} z_\nu z_\mu, 
\quad  z_\mu z_\nu^* = \lambda_{\nu\mu} z_\nu^* z_\mu,
\quad z_\mu^* z_\nu^* = \lambda_{\mu\nu} z_\nu^* z_\mu^*, \quad \mu,\nu = 0,1,2 ,
\ee
together with the spherical relation $\sum_\mu z_\mu^* z_\mu=1$. 
For $\theta=0$ one recovers the $*$-algebra of complex polynomial functions on the usual $S^4$.

The differential calculus $\Omega(\S^4)$ is generated as a graded differential $*$-algebra by the elements $z_\mu, z_\mu^*$ in degree 0 and elements $d z_\mu, d z_\mu^*$ in degree 1 with relations,
\be\label{rel:diff}
\begin{aligned}
&dz_\mu dz_\nu+ \lambda_{\mu\nu} dz_\nu dz_\mu =0, \\
&z_\mu dz_\nu = \lambda_{\mu\nu} dz_\nu z_\mu,
\end{aligned}
\qquad
\begin{aligned}
&dz_\mu dz_\nu^* + \lambda_{\nu\mu} dz_\nu^* dz_\mu =0;\\
&z_\mu dz_\nu^* = \lambda_{\nu\mu} dz_\nu^* z_\mu ,
\end{aligned}
\ee
and $\lambda_{\mu\nu}$ is as before. There is a unique differential $\dd$ on $\Omega(\S^4)$ such that  $\dd:z_\mu \mapsto d z_\mu$ and the involution on $\Omega(\S^4)$ is the graded extension of $z_\mu \mapsto z_\mu^*$:
$(\dd \omega)^*=\dd \omega^*$ and $(\omega_1\omega_2)^* = (-1)^{d_1 d_2}\omega_2^* 
\omega_1^*$ for $\omega_j$ a form of degree $d_j$.

The algebra $A(\Sk^7)$ of polynomial functions on the sphere $\Sk^7$ is generated by elements  
$\psi_a, \psi_a^*$, $a=1,\dots,4$, subject to relations
\be\label{s7t}
\psi_a \psi_b = \lambda'_{a b} \psi_b \psi_a, \quad  \psi_a \psi_b^* = \lambda'_{b a} \psi_b^* \psi_a,
\quad \psi_a^*\psi_b^* = \lambda'_{a b} \psi_b^* \psi_a^* ,
\ee
and with the spherical relation $\sum_a \psi_a^* \psi_a=1$. 
Now $\lambda'_{a b} = e^{2 \pi \ii \theta'_{ab}}$ and $(\theta'_{ab})$ is a real antisymmetric matrix; when it vanishes one gets the $*$-algebra of 
 complex polynomial functions on the sphere $S^7$. As before, a differential calculus $\Omega(\Sk^7)$ can be defined to be generated by the elements $\psi_a, \psi_a^*$ in degree 0 and elements $d \psi_a, d \psi_a^*$ in degree 1 satisfying relations similar to the ones 
in \eqref{rel:diff}.

A particular family of noncommutative 7 dimensional sphere $\Sk^7$
needs to be select for the principal noncommutative fibration over the given 4-sphere $\S^4$. We take the ones corresponding to the deformation parameters given by 
\be \label{def:lambda}
\lambda'_{ab}= 
\begin{pmatrix} 1 & 1 & \bar\mu & \mu \\ 
1 & 1 & \mu & \bar\mu \\
\mu &\bar\mu &1 & 1\\ 
\bar\mu & \mu &1 & 1 
\end{pmatrix}, \quad \mu = \sqrt{\lambda}, \qquad \mathrm{or} \qquad
\theta'_{ab}=\frac{\theta}{2}\begin{pmatrix} 0 & 0 & -1 & 1 \\ 
0 & 0 & 1 & -1 \\
1 & -1 & 0 & 0 \\ 
-1 & 1 & 0 & 0  \end{pmatrix}.
\ee
The previous choice is essentially the only one that  allows the algebra $A(\Sk^7)$ to carry an action 
 of the group $\SU(2)$ by automorphisms and such that the invariant subalgebra coincides with 
 $A(\S^4)$.  The best way to see this fact is by means of the matrix-valued function on $A(\Sk^7)$ given by 
\be\label{Psi} 
\Psi  =
\begin{pmatrix} 
\psi_1 & - \psi^*_2 \\ 
\psi_2 & \psi^*_1 \\
\psi_3 & -\psi^*_4 \\
\psi_4& \psi^*_3
\end{pmatrix}.
\ee
Then $\Psi^\dagger \Psi = \II_2 $ and 
\[
p = \Psi \Psi^\dagger
\] 
is a projection,  $p^2=p=p^\dagger$, with 
entries in $A(\S^4)$. Indeed, the right action of $\SU(2)$ on $A(\Sk^7)$ is simply given by
\be \label{actionSU2}
\alpha_w (\Psi) = \Psi w , \qquad w = \begin{pmatrix} w_1 & -\bar{w}_2 \\ w_2 & \bar{w}_1 \end{pmatrix} \in \SU(2), 
\ee
from which the invariance of the entries of $p$ follows at once. Explicitly,  
\be \label{projection1}
p= \half \begin{pmatrix}
1+z_0 & 0 & z_1 & - \bar{\mu} z_2^* \\
0 & 1+z_0 & z_2  & \mu z_1^* \\
z_1^*& z_2^* & 1-z_0 & 0\\
-\mu z_2 & \bar{\mu}   z_1 & 0 & 1-z_0 
\end{pmatrix}, 
\ee
with the generators of $A(\S^4)$ given in terms of the ones of $A(\Sk^7)$ by
\be\label{subalgebra}
\begin{aligned} 
z_0 &= \psi^*_1 \psi_1 + \psi^*_2 \psi_2 - \psi^*_3 \psi_3 - \psi^*_4 \psi_4&  \\
 &= 2(\psi^*_1 \psi_1 + \psi^*_2 \psi_2) -1 = 1 - 2(\psi^*_3 \psi_3 + \psi^*_4 \psi_4),   \\
z_1 &= 2 (\mu \psi_3^* \psi_1 + \psi^*_2 \psi_4)=2(\psi_1 \psi_3^* + \psi^*_2 \psi_4) ,  \\
z_2 &= 2(- \psi^*_1 \psi_4 +\bar\mu \psi_3^* \psi_2)=2(- \psi^*_1 \psi_4 + \psi_2 \psi_3^*).
\end{aligned}
\ee
One straightforwardly computes that $z_1^* z_1 + z_1^* z_1  + z_0 ^2 = 1$ and  the commutation rules 
$z_1 z_2  = \lambda z_2 z_1$, $z_1 z_2^* = \bar{\lambda} z_2^* z_1$, and that  $z_0$ is central.

The relations \eqref{subalgebra} can also be expressed in the form, 
\[ 
z_\mu =\sum_{ab}\psi^*_a(\gamma_\mu)_{ab}\psi_b,\qquad z_\mu^* = \sum_{ab}\psi^*_a(\gamma^*_\mu)_{ab}\psi_b,
\]
with $\gamma_\mu$ twisted $4 \times 4$ Dirac matrices given by
\begin{align}
\label{eq:dirac}
\gamma_0=\left(\begin{smallmatrix} 1 &&&\\ &1&&\\&&-1&\\&&&-1 \end{smallmatrix}\right),
\qquad
\gamma_1=2\begin{pmatrix} 0 & \begin{smallmatrix} 0 & 0 \\ 0 & 1 \end{smallmatrix} \\  \begin{smallmatrix} \mu & 0 \\ 0 & 0 \end{smallmatrix}& 0 \end{pmatrix}, \qquad
\gamma_2=2\begin{pmatrix} 0 & \begin{smallmatrix} 0 & -1 \\ 0 & 0 \end{smallmatrix} \\  \begin{smallmatrix} 0 & \bar\mu \\ 0 & 0 \end{smallmatrix}& 0 \end{pmatrix}.
\end{align}
These matrices satisfy  twisted Clifford algebra relations,
\[
\gamma_\mu \gamma_\nu +\lambda_{\mu\nu} \gamma_\nu\gamma_\mu=0,  \qquad
\gamma_\mu \gamma_\nu^* +\lambda_{\nu\mu} \gamma_\nu^*\gamma_\mu=4\delta_{\mu\nu}; \qquad (\mu,\nu =1,2).
\]
and $\gamma_0$ is the grading 
$$
\gamma_0 = -\frac{1}{4}[\gamma_1,\gamma_1^*][\gamma_2,\gamma_2^*].
$$
The manifolds $\S^4$ and $\Sk^7$ carry compatible toric actions. The torus $\bT^2$ acts on the algebra $A(\S^4)$ as 
\be\label{eq:act-S4}
\sigma_s(z_0, z_1, z_2) = (z_0, e^{2\pi i s_1} z_1, e^{2\pi i s_2} z_2), \quad s\in\bT^2 ,
\ee
and this action is lifted to a double cover action on $A(\Sk^7)$. The double cover map  $p: \tilde\bT^2 \to \bT^2$ is  
given explicitly by $p:(s_1,s_2) \mapsto (s_1+s_2,-s_1+s_2)$.
Then $\tilde\bT^2$ acts on the $\psi_a$'s as
\be \label{eq:lift-S7}
\tilde\sigma: \left(  \psi_1, \psi_2, \psi_3, \psi_4 \right) 
\mapsto \left(e^{2\pi i s_1}~\psi_1, ~e^{-2\pi i s_1}~\psi_2, ~e^{-2\pi i s_2}~\psi_3, ~e^{2\pi i s_2}~\psi_4 \right) 
\ee
Equation \eqref{subalgebra} shows that $\tilde\sigma$ is indeed a lifting to $\Sk^7$ of the action of $\bT^2$ on $\S^4$. This compatibility is built in the construction of the Hopf fibration $\Sk^7 \to \S^4$ as a deformation of the classical Hopf fibration $S^7 \to S^4$ with respect to an action of $\bT^2$, a fact that also dictated the form of the deformation parameter $\lambda'$ in \eqref{def:lambda}.
As we shall see later on, the previous double cover of tori comes from a spin cover $\Spin_\theta(5)$ of $\SO_\theta(5)$ deforming the usual action of $\Spin(5)$ on $S^7$ as a double cover of the action of $\SO(5)$ on $S^4$.
We refer to \cite{LS04} for the construction of the algebra inclusion $\A(\S^4) \into \A(\Sk^7)$ as a noncommutative principal bundle. 

\subsection{Associated bundles}\label{se:associated-modules}
 
Let $\rho$ be any representation of $\SU(2)$ on an $n$-dimensional vector space $V$. 
The corresponding equivariant maps are defined as the cotensor product
\[
\E:=\Cinf(\Sk^7) \boxtimes_\rho V:= \big\{\eta \in \Cinf(\Sk^7) \otimes V: (\alpha_w\otimes \id) (\eta) =(\id \otimes \rho(w)^{-1})(\eta)\big\},
\]
where $\alpha_w$ is the $\SU(2)$ action on $\Sk^7$ in \eqref{actionSU2}. The previous $\E$ is clearly a $\Cinf(\S^4)$ bimodule. 
We have proved in \cite{LS04} that $\E$ is also a finite projective module. It is worth stressing that the choice of a projection for a finite projective module requires the choice of one of the two (left or right) module structures. Similarly, the definition of a Hermitian structure requires the choice of the left or right module structure. In the following, we will always work with the right structure for the associated modules.
There is a natural right Hermitian structure on $\E$, defined in terms of the inner product of $V$ as
\[
\langle \eta, \eta' \rangle := \sum_i \bar \eta_i \eta'_i, 
\]
where we denoted $\eta=\sum_i\eta_i \otimes e_i$, and $\eta'=\sum_i \eta'_i \otimes e_i $ for a basis $\{e_i\}_{i=1}^n$ in $V$. One quickly checks that $\langle \eta, \eta' \rangle$ is an element in $\Cinf(\S^4)$, and that $\langle \mb , \mb \rangle$ satisfies all conditions of a right Hermitian structure. 

The bimodules given by the cotensor product $\Cinf(\Sk^7) \boxtimes_\rho V$ are of the type described in Section~\ref{se:toric-ncm} as the result of a quantization map. Indeed, the associated vector bundle $E=S^7 \times_\rho V$ on $S^4$ carries an action of $\tilde\bT^2$ induced from its action on $S^7$, which is obviously $\sigma$-equivariant. By the very definition of $\Cinf(\Sk^7)$ and of  $\Gamma(\S^4,E)$,  it follows that $\Cinf(\Sk^7) \boxtimes_\rho V \isom \Gamma(\S^4,E)$. 
Indeed, from the isomorphism, $\Gamma(S^4,E)\isom \Cinf(S^7)\boxtimes_\rho V$, the quantization map $L_{\theta'}$ of $\Cinf(S^7)$ -- when acting only on the first leg of the cotensor product -- establishes this isomorphism, 
\[
L_{\theta'} : \Cinf(S^7) \boxtimes_\rho V \to \Cinf(\Sk^7)\boxtimes_\rho V . 
\]
The above is well defined since the action of $\tilde\bT^2$ commutes with the action of $SU(2)$. 
Also, it is such that $L_{\theta'}(f \eta)=L_{\theta'}(f) \lt_{\theta'} L_{\theta'} (\eta)=L_{\theta}(f) \lt_{\theta'} L_{\theta'} (\eta)$ for $f \in \Cinf(S^4)$ and $\eta \in \Cinf(S^7) \boxtimes_\rho V$ -- due to the identity $L_{\theta'}=L_{\theta}$ on $\Cinf(S^4)\subset\Cinf(S^7)$.
A similar result holds for the right action $\rt_{\theta'}$.  

\bigskip

Given the right $\Cinf(\S^4)$-module $\E$, the {\it dual module}  is defined by 
\[
\E' := \big\{ \phi: \E \to \Cinf(\S^4) : \phi(\eta f ) = \phi(\eta) f, \, \, f\in\Cinf(\S^4) \big\},  
\]
and is naturally a left $\Cinf(\S^4)$-module.  
If $\E:=\Cinf(\Sk^7) \boxtimes_\rho V$ comes from the $\SU(2)$-representation $(V,\rho)$, by using the 
induces dual representation $\rho'$ on the dual vector space $V'$ given by
\[
\big( \rho'(w) v'\big) (v) := v'\big(\rho(w)^{-1} v\big); \qquad \forall v' \in V', v \in V ,
\]
we have that
\[
\E' \isom \Cinf(\Sk^7) \boxtimes_{\rho'} V' := \big\{ \phi \in \Cinf(\Sk^7) \otimes V': (\alpha_w \otimes \id) (\phi) =(\id \otimes \rho'(w)^{-1})(\phi)\big\} .
\]
Next, let $L(V)$ denote the space of linear maps on $V$, so that $L(V)=V\otimes V'$. The adjoint action of $\SU(2)$ on $L(V)$ is the tensor product representation $\ad:=\rho \otimes \rho'$ on $V \otimes V'$. We have the corresponding 
\[ 
\Cinf(\Sk^7) \boxtimes_{\ad} L(V) := \big\{ T \in \Cinf(\Sk^7) \otimes L(V): (\alpha_w \otimes \id) (T) =(\id \otimes \ad(w)^{-1})(T)\big\}  ,
\]
and write $T=T_{ij} \otimes e_{ij}$ with respect to the basis $\{ e_{ij} \}$ of $L(V)$ induced from the basis $\{e_i\}_{i=1}^n$ of $V$ and the dual basis $\{e'_i\}_{i=1}^n$ of $V'$. 

On the other hand, we have the endomorphism algebra 
\[
\End_{\Cinf(\S^4)}(\E):= \big\{ T: \E \to \E : T(\eta f ) = T (\eta) f, \, \, f\in\Cinf(\S^4) \big\}.  
\]
\begin{prop} \label{prop:end}
There is an isomorphism of algebras 
\[
\End_{\Cinf(\S^4)}(\E) \isom \Cinf(\Sk^7) \boxtimes_{\ad} L(V).
\] 
\end{prop}
\begin{proof} 
Recall that $\E \otimes_{\Cinf(\S^4)} \E'\subset \End_{\Cinf(\S^4)}(\E)$ densely (in the operator norm, cf. for instance \cite{GVF01}). We define a map from $\End_{\Cinf(\S^4)}(\E)$ to $\Cinf(\Sk^7) \boxtimes_\ad L(V)$ on this dense 
subalgebra by
$$
\eta \otimes \phi \mapsto \eta_i \phi_j \otimes e_{ij}, 
$$
with $\eta=\eta_i \ot e_i$ and $\phi=\phi_i \ot e'_i$.
On the other hand, $T\in\Cinf(\Sk^7) \boxtimes_\ad L(V)$ acts on $\eta \in \E$ by
$$
(T,\eta) \mapsto (T_{ij} \eta_j) \otimes e_i,
$$
which is clearly a right $\Cinf(\S^4)$-linear map with image in $\E$. Hence, the opposite inclusion 
$\Cinf(\Sk^7) \boxtimes_\ad L(V) \subset \End_{\Cinf(\S^4)}(\E)$.
\end{proof}
We see that the algebra of endomorphisms of $\E$ can be understood as the space of sections of the noncommutative vector bundle associated to the adjoint representation on $L(V)$ -- exactly as it happens in the classical case. This also allows an identification of skew-hermitian endomorphisms $\End^s(\E)$ -- which are defined in general in \eqref{def:skew} -- for the toric deformations at hand. 
\begin{coro}
There is an identification
 \[
 \End^s(\E) \isom \Cinf_\bR(\Sk^7) \boxtimes_\ad u(n),
 \]
  with $\Cinf_\bR(\Sk^7)$ denoting the algebra of self-adjoint elements in $\Cinf(\Sk^7)$ whereas $u(n)$ consists of skew-adjoint matrices in $M_n(\C) \isom L(V)$. 
\end{coro}
\begin{proof}
Note that the involution $T \mapsto T^*$ in $\End_{\Cinf(\S^4)}(\E)$ reads in components $T_{ij} \mapsto \bar{T_{ji}}$ so that, with the identification of Proposition~\ref{prop:end}, the algebra $\End^s(\E)$ is made of elements $X \in \Cinf(\Sk^7) \boxtimes_\ad L(V)$ satisfying $\bar{X_{ji}} = -X_{ij}$. Since any element in $\Cinf(\Sk^7)$ can be decomposed as the sum of two self-adjoint elements, $X_{ij} = X_{ij}^\Re + \ii X_{ij}^\Im$, we can write 
$$
X=\sum_i X_{ii}^\Im\otimes \ii e_{ii}+ \sum_{i \neq j} X_{ij}^\Re \otimes (e_{ij} - e_{ji})  + X_{ij}^\Im \otimes (\ii e_{ij}+\ii e_{ji}) = \sum_a X_a \otimes \sigma^a,
$$
where $X_a$ are generic elements in $\Cinf_\bR(\Sk^7)$ and $\{\sigma^a, a=1,\ldots,n^2\} $ are the generators of $u(n)$. 
\end{proof}

\subsection{Yang--Mills theory on $\S^4$}
\label{se:ym}
Let us now discuss the Yang--Mills action functional on $\S^4$ together with its equations of motion. We will see that instantons naturally arise as absolute minima of this action. Before we proceed we recall the noncommutative spin structure $(\Cinf(\S^4),\ch,D,\gamma_5$) of $\S^4$ with $\ch=L^2(S^4,\cs)$ the Hilbert space of spinors, $D$ the undeformed Dirac operator, and $\gamma_5$ -- the even structure -- the fifth Dirac matrix.

Let $\E=\Gamma(\S^4,E)$ for some $\sigma$-equivariant vector bundle $E$ on $S^4$, so that there exists a projection $p \in M_N(\Cinf(\S^4))$ such that $\E\isom p (\Cinf(\S^4)^N$. Recall from Appendix~\ref{se:connections} that a connection $\nabla$ on $\E=\Gamma(\S^4,E)$  is a map from $\E$ to $\E\otimes \Omega(\S^4)$ satisfying a Leibniz rule. The Yang--Mills action functional is defined in terms of the curvature of the connection $\nabla$ on $\E$, which, as seen in the appendix, is an element in 
$\Hom_{\Cinf(\S^4)} (\E, \E \otimes \Omega^2(\S^4) )$. Equivalently, it is an element in $\End_{\Omega(\S^4)}(\E\ot \Omega(\S^4))$ of degree 2. An inner product on the latter algebra is defined as follows \cite[III.3]{C94}. Any element $T \in \End_{\Omega(\S^4)}(\E\ot \Omega(\S^4))$ of degree $k$ is understood as an element in $p M_N(\Omega^k(\S^4))p$, since $\E\ot \Omega(\S^4)$ is a finite projective module over $\Omega(\S^4)$. A trace over internal indices together with the inner product in \eqref{eq:inner-product-forms}, defines the inner product $(\cdot,\cdot)_2$ on $\End_{\Omega(\S^4)}(\E\ot \Omega(\S^4))$. Having that, we can make the following definition.
\begin{defi}
\label{def:YM}
The \emph{Yang--Mills action functional} on the collection $C(\ce)$ of compatible connections on $\E$ is defined by
$$
\YM(\nabla)=\big( F,F\big)_2 = \ncint *_\theta \tr (F *_\theta F),
$$
for any connection $\nabla$ with curvature $F$. 
\end{defi}
\noindent
Recall from Appendix~\ref{se:connections} that unitary endomorphisms of $\E$ make up gauge transformations for connections. 

\begin{lma} 
The Yang--Mills functional is gauge invariant, positive and quartic.
\end{lma}
\begin{proof} 
From equation~\eqref{ugcur}) the curvature $F$ transforms as $F \mapsto u^* F u$ under a gauge transformation $u \in \cU(\E)$. Since $\cU(\E)$ can be identified with the unitary elements in $p M_N(\A) p$, it follows that 
\begin{align*}
\YM(\nabla^u)=\ncint \sum_{i,j,k,l} *_\theta (\bar{u_{ji}} F_{jk} *_\theta 
F_{kl} u_{li} )=\YM(\nabla), 
\end{align*}
using the tracial property of the Dixmier trace and the fact that $u_{li} \bar{u_{ji}}=\delta_{lj}$. 
Positiveness follows from the equality 
$$
(F,F)_2=(F_D,F_D)_D=\ncint F_D^* F_D,
$$
and the last expression is clearly positive.
\end{proof}

A variational principle yields Yang--Mills equations. For this, it is enough to  consider linear perturbations $\nabla_t=\nabla + t \alpha$ of a connection $\nabla$ on $\E$ by elements $\alpha \in \Hom_{\Cinf(\S^4)}(\E,\E\ot_{\Cinf(\S^4)} \Omega^1(\S^4)) $. The curvature $F_t$ of $\nabla_t$ is readily computed as $F_t = F + t [\nabla, \alpha] + \cO (t^2)$. If we suppose that $\nabla$ is an extremum of the Yang--Mills action functional, this linear perturbation should not affect the action to first order. In other words, we should impose
\[
\frac{\partial}{\partial t}\bigg|_{t=0} \YM(\nabla_t)=0.
\]
If we substitute the explicit formula for $F_t$, we obtain
\[
\big( [\nabla, \alpha],F \big)_2 + \bar{\big( [\nabla, \alpha] , F\big)_2} = 0.
\]
Positive definiteness of the scalar product implies that 
$(F_t,F_t)_2=\bar{(F_t,F_t)_2}$, which differentiated with respect to $t$, at $t=0$, yields $\big( [\nabla, \alpha],F \big)_2 = \bar{\big( [\nabla, \alpha] , F\big)_2}$; hence, $( [\nabla, \alpha],F )_2=0$. Being $\alpha$ arbitrary, we derive the equations of motion
\[
[\nabla^*, F\big] = 0,
\]
with the adjoint of $[\nabla,\cdot]$ defined with respect to the inner 
product $(\cdot, \cdot)_2$, {i.e.}
\[
\big( [\nabla^*,\alpha], \beta \big)_2 = \big( \alpha, [\nabla, \beta] \big)_2 ,
\]
for $\alpha \in \Hom_{\Cinf(\S^4)}(\E,\E \ot \Omega^3(\S^4))$ and $\beta \in \Hom_{\Cinf(\S^4)}(\E,\E\ot \Omega^1(\S^4))$.
From Lemma~\ref{lma:d-dstar}, it follows that $[\nabla^* ,F]= \ast_\theta [\nabla, \ast_\theta F]$, so that the equations of motion can also be written as the more familiar \textit{Yang--Mills equations}:
\be \label{eq:ym}
[\nabla, \ast_\theta F]=0.
\ee
Note that connections with a self-dual or antiself-dual curvature $\ast_\theta F= \pm F$ are special solutions of the Yang--Mills equations. Indeed, in this case the latter equations follows directly from the Bianchi identity $[\nabla,F]=0$ in Proposition \ref{ubianchi}. 

\bigskip

Suppose $\E$ is a finite projective module over $\Cinf(\S^4)$ defined by a projection $p \in M_N(\Cinf(\S^4))$. The topological action for $\E$ is given by a pairing between the class of $p$ in K-theory and the cyclic cohomology of $\Cinf(\S^4)$ \cite[VI.3]{C94}. For computational purposes, we give a  definition of the topological action in terms of the curvature of a connection on $\E$
\begin{defi}\label{topact}
Let $\nabla$ be a connection on $\E$ with curvature $F$. The {\rm topological action} is given by 
$$
\Top(\E)= (F,\ast_\theta F)_2= \ncint \ast_\theta \tr(F^2)
$$
where the trace is taken over internal indices and in the second equality we have used the identity $\ast_\theta \circ \ast_\theta = \id$ on $\S^4$.
\end{defi}
This definition does not depend on the choice of a connection on $\E$. Since two connections differ by an element $\alpha \in \Hom_{\Cinf(\S^4)}(\E,\E\otimes \Omega^1(\S^4))$, we need to establish that $(F',\ast_\theta F')_2=(F,\ast_\theta F)_2$ where $F'=F+t [\nabla,\alpha] + \cO(t^2)$ is the curvature of $\nabla':=\nabla + t \alpha$, $t \in \bR$. By the definition of the inner product it follows
\begin{align*}
(F', \ast_\theta F')_2 - (F, \ast_\theta F)_2 &= t (F, \ast_\theta [\nabla, \alpha])_2+ t([\nabla,\alpha],\ast_\theta F)_2 +\cO(t^2) \nn \\
&= t(F, [\nabla^*,*_\theta \alpha])_2+ t([\nabla^*,*_\theta\alpha],F)_2 +\cO(t^2),
\end{align*}
which vanishes due to the Bianchi identity $[\nabla,F]=0$.

The Hodge star operator $\ast_\theta$ splits $\Omega^2(\S^4)$ into a self-dual and an antiself-dual parts,
\[
\Omega^2(\S^4)=\Omega_+^2(\S^4) \op \Omega_-^2(\S^4).
\]
In fact, $\Omega_\pm^2(\S^4)=L_\theta \left(\Omega_\pm^2(S^4) \right)$. This decomposition is orthogonal with respect to the inner product $( \cdot,\cdot)_2$, so that we can write the Yang--Mills action functional as
\[
\YM(\nabla) =  \big( F_+,F_+\big)_2 + \big( F_-,F_-\big)_2.
\]
Comparing this with the topological action,
\[ 
\Top(\E)= \big( F_+,F_+\big)_2 - \big( F_-,F_-\big)_2,
\]
we see that $\YM(\nabla) \geq | \Top(\E) |$, with equality holding iff 
\[
\ast_\theta F = \pm F. 
\]
Solutions of these equations are called instantons. We conclude that instantons are absolute minima of the Yang--Mills action functional.

\section{$\SU(2)$-instantons on $\S^4$}\label{sect:constrinst}

In this section, we describe the  family of $\SU(2)$ instantons of charge $1$ on $\S^4$, constructed in \cite{LS06} by acting with a twisted infinitesimal conformal symmetry on a basic instanton on $\S^4$. 
The twisting is implemented with a twist of Drinfel'd type \cite{Dri90} -- in fact, explicitly constructed by Reshetikhin \cite{Res90}. 
There it was also proven that the `tangent space' of the moduli space of irreducible instantons at the basic instanton is five-dimensional, thus this set is complete. 
We perturb the connection $\nabla_0$ of the basic instanton found in \cite{CL01},  linearly by sending $\nabla_0 \mapsto \nabla_0+ t \alpha$ where $\alpha \in \Hom_{\Cinf(\S^4)}(\E,\E \ot \Omega^1(\S^4))$, and $t\in \bR$. For this new connection to still be an instanton, we have to impose the self-duality equations on its curvature. By deriving with respect to $t$, at $t=0$, we obtain the linearized self-duality equations to be fulfilled by $\alpha$. It is in this sense that we are considering the tangent space to the moduli space of instantons at $\nabla_0$. 

Thus,  of central interest is the noncommutative instanton bundle first constructed in \cite{CL01}.
Now $V=\C^2$ and $\rho$ is the defining representation of $\SU(2)$. 
The projection $p$ giving the $\Cinf(\S^4)$-module $\Cinf(\Sk^7)\boxtimes_\rho \C^2$ as a direct summand of $\big( \Cinf(\S^4) \big)^N$ for some $N$, is precisely given by the one in \eqref{projection1} and $N=4$. Indeed, a generic element in $\Cinf(\Sk^7)\boxtimes_\rho \C^2$ is of the form $\Psi^\dagger f$ for some $f \in \Cinf(\S^4) \otimes \C^4$ with $\Psi$ defined in \eqref{Psi}, and the correspondence is given by 
\[ 
\Cinf(\Sk^7)\boxtimes_\rho \C^2  \isom p \big( \Cinf(\S^4) \big)^4, \quad \Psi^\dagger f \leftrightarrow p f.
\]
Furthermore, $\End_{\Cinf(\S^4)}(\E) \isom \Cinf(\Sk^7)\boxtimes_{\ad} M_2(\C))$. It is a known fact that $M_2(\C)$ decomposes into the adjoint representation $su(2)$ and the trivial representation $\C$ while it is easy to see that  $\Cinf(\Sk^7 )\boxtimes_{\id} \C \isom \Cinf(\S^4)$. Thus, we conclude that 
\[
\End_{\Cinf(\S^4)}(\E) \isom \Gamma(\ad(\Sk^7)) \oplus \Cinf(\S^4),
\]
where we have set $\Gamma(\ad(\Sk^7)):=\Cinf(\Sk^7)\boxtimes_{\ad} su(2)$. The latter $\Cinf(\S^4)$-bimodule will be understood as the space of (complex) sections of the adjoint bundle. It is the complexification of the traceless skew-hermitian endomorphisms $\Cinf_\bR(\Sk^7) \boxtimes_\ad su(2)$.

We let $\nabla_0=p\circ \dd$ be the canonical (Grassmann) connection on the projective module 
$\E=p (\Cinf(\S^4) )^4 \simeq \Cinf(\Sk^7)\boxtimes_\rho \C^2$, with the projection $p= \Psi^\dagger \Psi$ of \eqref{projection1} and $\Psi$ is the matrix \eqref{Psi}. On equivariant maps we write the connection $\nabla_0$ as
\be\label{cancon}
\nabla_0: \E \to \E \ot_{\Cinf(\S^4)} \Omega^1(\S^4),  \qquad
(\nabla_0 f)_i = \dd f_i + \omega_{ij} \times_\theta f_j,
\ee
where $\omega$ -- referred to as the gauge potential -- is given in terms of the matrix $\Psi$ by
\be\label{cangp}
\omega=\Psi^\dagger \dd \Psi . 
\ee
The above, is a $2\times 2$-matrix with entries in $\Omega^1(\Sk^7)$ satisfying $\bar{\omega_{ij} }=\omega_{ji}$ and $\sum_i \omega_{ii}=0$. 
Note here that the entries $\omega_{ij}$ commute with all elements in $\Cinf(\Sk^7)$. Indeed, from \eqref{Psi} we see that the elements in $\omega_{ij}$ are $\bT^2$-invariant and hence central (as one forms) in $\Omega(\Sk^7)$. 
In other words $L_\theta(\omega)=\omega$, showing that for an element $f \in \E$ as above, we have $\nabla_0(f)_i=\dd f_i + \omega_{ij} \times_\theta f_j= \dd f_i + \omega_{ij} f_j$ which coincides with the action of the classical connection $\dd+\omega$ on $f$.

The curvature $F_0=\nabla_0^2=\dd \omega + \omega^2$ of $\nabla_0$  satisfies the self-duality equations  \cite{AB02,CD02},
\[
\ast_\theta F_0 = F_0;
\] 
hence this connection is an instanton. At $\theta=0$, the connection \eqref{cangp} is nothing but the $SU(2)$ instanton of \cite{BPST75}. 
Its `topological charge', i.e. the values of $\Top(\E)$ in Definition \ref{topact}, was already computed in \cite{CL01}. It depends only on the class $[p]$ of the bundle and can be evaluated as the index 
\be\label{topbasic}
\Top([p]) = \ind (D_p) = \ncint \gamma_5 \pi_D(\chern_2(p))
\ee
 of the twisted Dirac operator 
\[
D_p = p (D \ot \II_4) p .
\]
The last equality in \eqref{topbasic} follows from the vanishing of the first Chern character class $\chern_1(p)$ of the bundle (the Chern character classes and their realization as operators are in Appendix~\ref{se:cc}). 
On the other hand, for the second class one finds 
\[
\pi_D \big(\chern_2(p))\big) =  3 \gamma_5, 
\]
which, together with the fact that 
\[
\ncint 1 = \dix |D|^{-4} = \frac{1}{3},
\] 
on $S^4$ (see for instance \cite{GVF01,Lnd97}),
gives  the value $\Top([p])=1$. 

\bigskip
We aim at constructing all connections $\nabla$ on $\E$ whose curvature satisfies the self-duality equations and of topological charge equal to $1$. 
We can write any such connection in terms of the canonical connection as in equation~\eqref{uconn}, {i.e.} $\nabla=\nabla_0 + \alpha$ with $\alpha$ a one-form valued endomorphism of $\E$. Clearly, this will not change the value of the topological charge. We are particularly interested in $\SU(2)$-instantons, so we impose that $\alpha$ is traceless and skew-hermitian. Here the trace is taken in the second leg of $\End_{\Cinf(\S^4)}(\E) \isom \Sk^7) \boxtimes_\ad M_2(\C)$. When complexified, this gives an element $\alpha \in \Omega^1(\S^4)\ot_{\Cinf(\S^4)}\Gamma(\ad(\Sk^7)) =: \Omega^1(\ad(\S^4))$.

\bigskip
The noncommutative sphere $\S^4$ carries \cite{Sit01} a twisted symmetry action of the enveloping algebra $\U_\theta(so(5))$. This is lifted to $\Sk^7$ and the above basic instanton $\nabla_0$ is invariant under this infinitesimal quantum symmetry. 

Different instantons are obtained by a twisted symmetry action of $\U_\theta(so(5,1))$. Classically, $so(5,1)$ is the conformal Lie algebra consisting of the infinitesimal diffeomorphisms leaving the conformal structure invariant. We construct the Hopf algebra $\U_\theta(so(5,1))$ by adding 5 generators to $\U_\theta(so(5))$ and describe its action on $\S^4$ together with its lift to $\Sk^7$. The induced action of $\U_\theta(so(5,1))$ on forms leaves the conformal structure invariant. The action of  $\U_\theta(so(5,1))$ on $\nabla_0$  eventually results in a five-parameter family of instantons. 

\subsection{Twisted infinitesimal  symmetries}
We start with the construction of the twisted symmetry $\U_\theta(so(5))$. The eight roots of the Lie algebra $so(5)$ are two-component vectors $r=(r_1,r_2)$ 
given explicitly by
$r=(\pm1,\pm1)$ and $r=(0,\pm1)$, $r=(\pm1,0)$.  There are corresponding generators $E_r$ of $so(5)$ together 
with two generators $H_1,H_2$ of the abelian Cartan subalgebra. The Lie brackets are
\be\label{lie-so5}
\begin{aligned}
& [H_1,H_2] = 0, \quad [H_j,E_r] = r_j E_r , \\ 
&  [E_{-r},E_{r}] = r_1 H_1 + r_2 H_2, \quad 
[E_{r},E_{r'}] = N_{r,r'} E_{r+r'}, 
\end{aligned}
\ee
with $N_{r,r'}=0$ if $r+r'$ is not a root.
The universal enveloping algebra $\U(so(5))$ is the algebra 
 generated by elements $\{H_j, E_r\}$ modulo relations given by the previous Lie brackets. It is a Hopf algebra with, on the generators, the coproduct 
\[ 
\Delta_0 : \U(so(5)) \to \U(so(5)) \ot \U(so(5)), \quad 
\quad X \mapsto \Delta_0(X)=X \ot \II + \II \ot X , 
\]
the counit $\varepsilon(X)=0$ and the antipode $S(X)=-X$ .
 
The twisted universal enveloping algebra $\U_\theta(so(5))$ is 
generated as above (i.e. we do not change the algebra structure) but is endowed with a twisted 
coproduct, 
\[
\Delta_\theta: \U_\theta(so(5)) \to \U_\theta(so(5)) \ot \U_\theta(so(5)), \quad
\quad X \mapsto \Delta_\theta(X) = \cf \Delta_0(X) \cf^{-1}. 
\]
For the symmetries of the present paper the twist $\cf$ is taken as 
\[
\cf = \lambda^{-H_1 \ot H_2} .
\]
On the generators $E_r$, $H_j$, the twisted coproduct reads
\be\label{twdel} 
\begin{aligned}
& \Delta_\theta(E_r) = E_r \ot \lambda^{-r_1 H_2} 
+ \lambda^{-r_2 H_1} \ot E_r , \\
& \Delta_\theta(H_j) = H_j  \ot \II + \II \ot H_j .
\end{aligned}
\ee
This coproduct allows one to represent $\U_\theta(so(5))$ as an algebra of twisted derivations on both $\S^4$ and 
$\Sk^7$ as we shall see  below. With counit and antipode given
\be\label{twhopf} 
\begin{aligned}
& \varepsilon(E_r) = \varepsilon(H_j) = 0, \\
& S(E_r) = - \lambda^{r_2 H_1}  E_ r\lambda^{r_1 H_2}, \quad S(H_j) = -H_j ,
\end{aligned}
\ee
 the algebra $\U_\theta(so(5))$ becomes a Hopf algebra (see for instance \cite{CP94}). At the classical value of the deformation 
parameter, $\theta=0$, one recovers the Hopf algebra structure of $\U(so(5))$.

 Let us introduce `partial derivatives', 
$\partial_\mu$ and 
 $\partial_\mu^*$ with the usual action on the generators of the algebra $\A(\S^4)$ {i.e},  
 $\partial_\mu(z_\nu)=\delta_{\mu\nu}$, $\partial_\mu(z_\nu^*)=0$, and 
 $\partial_\mu^*(z_\nu^*)=\delta_{\mu\nu}$, $\partial_\mu^*(z_\nu)=0$. Then, the action of $\U_\theta(so(5))$ on 
$\A(\S^4)$ is given by the operators,
\be\label{act4}
\begin{aligned}
H_1 &= z_1 \partial_1 - z_1^* \partial_1^* , \\
E_{+1,+1}&= z_2 \partial_1^*   - z_1 \partial_2^*,\\
E_{+1,0} &= \tfrac{1}{\sqrt{2}} (2 z_0 \partial_1^* - z_1 \partial_0), 
\end{aligned}
\qquad
\begin{aligned}
H_2 &= z_2 \partial_2 - z_2^* \partial_2^* \\
E_{+1,-1} &= z_2^* \partial_1^* - z_1 \partial_2  \, , \\
E_{0,+1} &=  \tfrac{1}{\sqrt{2}} (2 z_0 \partial_2^* - z_2 \partial_0) \, , 
\end{aligned}
\ee
and $E_{-r}=(E_{r})^*$, with the obvious meaning of the adjoint.  
These operators are extended to the whole of $\A(\S^4)$ as twisted derivations via the coproduct \eqref{twdel}:
\be\label{tder} 
\begin{aligned}
& E_r( a b ) := m \Big( \Delta_\theta(E_r) (a\ot b) \Big) 
= E_r(a) \lambda^{-r_1 H_2}(b) + \lambda^{-r_2 H_1}(a) E_r(b) , \\
& H_j( a b ) := m \Big( \Delta_\theta(H_j) (a\ot b) \Big) = H_j(a) b + a H_j(b) , 
\end{aligned}
\ee
for any two elements $a,b\in \A(\S^4)$, and $m(a\ot b):=ab$ is the algebra multiplication.  With these twisted rules, one readily 
checks compatibility with the commutation relations \eqref{s4t} of $\A(\S^4)$. 

The representation of $\U_\theta(so(5))$ on $\S^4$ given in \eqref{act4} is the fundamental vector representation.  
 When lifted to $\Sk^7$ one gets the fundamental spinor representation: as we see from the 
 quadratic relations among corresponding generators, given in \eqref{subalgebra}, the lifting amounts to 
take the `square root' representation.
 The action of $\U_\theta(so(5))$ on $\A(\Sk^7)$ is constructed by requiring twisted derivation properties via the 
 coproduct \eqref{tder} so as to reduce to the action \eqref{act4} on $\A(\S^4)$ when using the defining 
 quadratic relations \eqref{subalgebra}. The action on  $\A(\Sk^7)$ can be given as the action of matrices 
$\Gamma$'s on the $\psi$'s:
\begin{align}\label{act7}
\psi_a \mapsto \sum_b \Gamma_{ab} \psi_b; \qquad \psi^*_a \mapsto  \sum_b \tilde\Gamma_{ab} \psi^*_b,
\end{align}
with the matrices $\Gamma = \{H_j, E_r\}$ given explicitly by,
\be\label{tgamma}
\begin{aligned}
&H_1 = \half\left( \begin{smallmatrix}
1  & & & \\
 & -1 &  &   \\
 &   & -1 &  \\
 &  & & 1 
\end{smallmatrix} \right), \\
&E_{+1,+1} =\begin{pmatrix} 
0 & 0 \\ 
0 & \begin{smallmatrix} 0 & -1 \\ 0 & 0 \end{smallmatrix} \\
\end{pmatrix},  \\
&E_{+1,0} =  \tfrac{1}{\sqrt{2}} \begin{pmatrix} 
0 & \begin{smallmatrix} 0 & 0\\ 0 & -1 \end{smallmatrix} \\ 
\begin{smallmatrix} \mu  & 0 \\ 0 & 0 \end{smallmatrix} & 0 \\
\end{pmatrix}, 
\end{aligned}
\qquad
\begin{aligned}
&H_2 = \half\left( \begin{smallmatrix}
-1  & & & \\
 & 1 &  &   \\
 &   & -1 &  \\
 &  & & 1 
\end{smallmatrix} \right), \\
&E_{+1,-1} =\begin{pmatrix} 
\begin{smallmatrix} 0 & 0 \\ -\mu & 0 \end{smallmatrix} & 0 \\ 
0 & 0 \\
\end{pmatrix},  \\
&E_{0,+1} =  \tfrac{1}{\sqrt{2}} \begin{pmatrix} 
0 & \begin{smallmatrix} 0 & \bar{\mu} \\ 0 & 0 \end{smallmatrix} \\ 
\begin{smallmatrix} 0  & 1 \\ 0 & 0 \end{smallmatrix} & 0 \\
\end{pmatrix}, 
\end{aligned}
\ee
and $\tilde\Gamma:= \sigma \Gamma \sigma^{-1}$ with
\[
\sigma := \begin{pmatrix} \begin{smallmatrix} 0 & -1 \\ 1 & 0 \end{smallmatrix} & 0 \\  0 & \begin{smallmatrix} 0 & -1 \\ 1 & 0 \end{smallmatrix} \end{pmatrix}.
\]
Furthermore, $E_{-r}=(E_{r})^*$. Notice that $\tilde\Gamma=-\Gamma^t$ at $\theta=0$. 
With the twisted rules \eqref{tder} for the action on products, one checks compatibility of the above 
action with the commutation relations \eqref{s7t} of $\A(\Sk^7)$. 
\begin{rema}
Note that the operators $H_1$ and $H_2$ in \eqref{act4} are the infinitesimal generators of the action of $\bT^2$ on $\S^4$ as given in equation~\eqref{eq:act-S4}.
On the other hand, by comparing the form of the matrices $H_1$ and $H_2$ in \eqref{tgamma} with the lifted action $\tilde \sigma$ of $\tilde \bT^2$ on $\Sk^7$ as defined in \eqref{eq:lift-S7}, one checks that 
\[
\tilde\sigma_s = e^{ \pi i \left((s_1+s_2) H_1 + (-s_1+s_2) H_2\right)} ,
\] 
when acting on the spinor $(\psi_1, \cdots, \psi_4)^t$.
\end{rema}
The relation between the spinor representation in \eqref{tgamma} and the twisted Dirac matrices in 
\eqref{eq:dirac} is easily established,
\[
\begin{aligned}
&\tfrac{1}{4}[\gamma_1^*, \gamma_1]=2 H_1\\
&\tfrac{1}{4}[\gamma_1, \gamma_2]=(\mu+\bar\mu) E_{+1,+1}\\
&\tfrac{1}{4}[\gamma_1, \gamma_0]= \sqrt{2}  E_{+1,0}
\end{aligned}
\qquad
\begin{aligned}
&\tfrac{1}{4}[\gamma_2^*, \gamma_2]=2 H_2\\
&\tfrac{1}{4}[\gamma_1, \gamma_2^*]=(\mu+\bar\mu) E_{+1,-1}\\
&\tfrac{1}{4}[\gamma_2, \gamma_0]=\sqrt{2}\bar\mu E_{0,+1}
\end{aligned}
\]
\begin{prop}
The instanton gauge potential $\omega$ is invariant under the twisted action of $\U_\theta(so(5))$. 
\end{prop}
\begin{proof}
One finds that the gauge potential transforms as
$$
\omega = \Psi^* \dd \Psi \mapsto \Psi^* \big( \tilde\Gamma^t \lambda^{-r_1 H_2} + \lambda^{r_2 H_1} \Gamma \big) \dd \Psi.
$$
with  
$\lambda^{-r_i H_j}$ understood in its representation \eqref{tgamma} on $\A(\Sk^7)$. Direct computation for $\Gamma=\{H_j, E_r\}$ shows that $\tilde\Gamma^t \lambda^{-r_1 H_2} + \lambda^{r_2 H_1} \Gamma=0$, which finishes the proof.
\end{proof}

\subsection{Twisted conformal transformations}\label{subsect:infconf}
In order to have new instantonic configurations we need to use conformal transformations. 
In the definition of the enveloping algebra $\U_\theta(so(5,1))$ we do not change the algebra structure, {i.e.} we take the relations of $\U(so(5,1))$, as we did in the case of $\U_\theta(so(5))$. We thus define $\U_\theta(so(5,1))$ as the algebra $\U_\theta(so(5))$ with five extra generators adjoined,  $H_0, G_r$, 
$r=(\pm 1,0),(0,\pm 1)$, subject to the relations of $\U_\theta(so(5))$ of equation~\eqref{lie-so5} together with the (undeformed) relations,
\begin{align*} 
\begin{aligned}
&[H_0,H_i]=0, \\ & [H_0,G_r]=\sqrt{2} E_r, 
\end{aligned} \qquad 
\begin{aligned}
&[H_j, G_r] = r_j G_r, \\ & [H_0,E_r]=(\sqrt{2})^{-1} G_r,
\end{aligned}
\end{align*}
whenever $r=(\pm 1,0),(0,\pm 1)$, and 
\begin{align*}
\begin{aligned}
&[G_{-r},G_r]=2 r_1 H_1 + 2 r_2 H_2, \\ & [E_r,G_{r'}]=\tilde N_{r,r'} G_{r+r'}, 
\end{aligned}
\qquad
\begin{aligned}
&[G_r,G_{r'}] =N_{r,r'} E_{r+r'}, \\ &[E_{-r},G_r]=\sqrt{2} H_0,
\end{aligned}
\end{align*}
with  as before, $N_{r,r'}=0$ if $r+r'$ is not a root of $so(5)$ and also $\tilde N_{r,r'}=0$ if 
$r+r' \notin \{(\pm 1,0),(0,\pm 1) \}$.
Although the algebra structure is unchanged, again the Hopf algebra structure of gets twisted. The twisted structures are given by equations \eqref{twdel} and \eqref{twhopf} together with analogous ones for the extra generators,
\[
\begin{aligned}
\Delta_\theta(G_r) &= G_r \ot \lambda^{-r_1 H_2} + \lambda^{-r_2 H_1} \ot G_r , \\ 
S(G_r) &= - \lambda^{r_2 H_1}  G_ r\lambda^{r_1 H_2}, \\
\varepsilon(G_r) &= 0 ,
\end{aligned}
\qquad
\begin{aligned}
\Delta_\theta(H_0) &= H_0 \ot 1+ 1 \ot H_0 ,\\
S(H_0) &=-H_0, \\
\varepsilon(H_0) &= 0
\end{aligned}
\]
As for $\U_\theta(so(5))$, these structures make $\U_\theta(so(5,1))$ a Hopf algebra. 

The action of $\U_\theta(so(5,1))$ on $\A(\S^4)$ is given by the operators \eqref{act4} and 
\begin{align}\label{confact4}
&H_0 = \partial_0 - z_0 (z_0 \partial_0 + z_1 \partial_1 + z_1^* \partial_1^* +
 z_2 \partial_2 + z_2^* \partial_2^*), \nn \\
&G_{1,0}=  2 \partial_1^* - z_1 (z_0 \partial_0 + z_1 \partial_1 + z_1^* \partial_1^* + \bar\lambda z_2 \partial_2 + \lambda z_2^* \partial_2^*),\\
&G_{0,1}=  2 \partial_2^* - z_2 (z_0 \partial_0 + z_1 \partial_1 + z_1^* \partial_1^* + z_2 \partial_2 + z_2^* \partial_2^*), \nn
\end{align}
with $G_{-r} = (G_r)^*$. The extra $\lambda$'s in $G_{1,0}$ (and $G_{-1,0}$) are necessary to preserve the Lie algebra structure of $\U_\theta(so(5,1))$ as can be seen by a direct computation. 

Again, the operators $H_0, G_r$ in\eqref{confact4} are extended to the whole of $\A(\S^4)$ by analogues of \eqref{tder} using the twisted coproduct, 
\bean
&& G_r( a b ) := m \Big( \Delta_\theta G_r (a\ot b) \Big) = G_r(a) \lambda^{-r_1 H_2}(b) + \lambda^{-r_2 H_1}(a) G_r(b) , \nn\\
&& H_0( a b ) := m \Big( \Delta_\theta H_0 (a\ot b) \Big) = H_0(a) b + a H_0(b), \qquad \forall \quad a,b \in \A(\S^4).
\eean

The Hopf algebra $\U_\theta(so(5,1))$ consists of the infinitesimal twisted conformal transformations on $\S^4$. Firstly, one extends 
the twisted action of $\U_\theta(so(5,1))$ on $\A(\S^4)$ to the differential calculus $(\Omega(\S^4),\dd)$ by requiring it to commute with the exterior derivative, 
\[
T \cdot \dd \omega := \dd (T \cdot \omega) .
\]
for $T\in  \U_\theta(so(5,1)), ~ \omega \in \Omega(\S^4)$. Then one has the 
\begin{lma}
\label{lma:so51-hodge}
With the above twisted action, the Hopf algebra $\U_\theta(so(5,1))$ leaves the Hodge $\ast_\theta$-structure of $\Omega(\S^4)$ invariant: 
$$
T\cdot (\ast_\theta\omega)= \ast_\theta (T \cdot \omega) .
$$
\end{lma}
\begin{proof}
This follows from the fact that $T(L_\theta(a))=L_\theta(t\cdot a)$ for $a \in \A(S^4)$ and $t \in \U(so(5,1))$ the classical limit ($\theta =0$) of $T \in \U_\theta(so(5,1))$. Then, since $\U(so(5,1))$ leaves the Hodge $\ast$-structure of $\Omega(S^4)$ invariant
and the differential $\dd$ commutes with the action of $\U_\theta(so(5,1))$, if follows that the latter algebra leaves the Hodge $\ast_\theta$-structure of $\Omega(\S^4)$ invariant as well.
\end{proof}
In the same manner as for $\U_\theta(so(5))$, the action of the Hopf algebra $\U_\theta(so(5,1))$ on $\S^4$ can be lifted to an action on $\Sk^7$. 
Again, the latter action can be written as in \eqref{act7} in terms of matrices $\Gamma$'s acting on the $\psi$'s,
where together with the matrices \eqref{tgamma} 
we have the additional matrices $\Gamma=\{H_0,G_r\}$,
\begin{align*}
H_0 &= \half (-z_0 \I_4 + \gamma_0), \nn\\
G_{1,0}&=\half (-z_1 \lambda^{-H_2} + \gamma_1), \qquad 
G_{0,1}=\half (-z_2 + \lambda^{-H_1} \gamma_2),  
\end{align*}
with $G_{-r}= (G_r)^*$ and $\tilde\Gamma = \sigma \Gamma \sigma^{-1}$. 
Notice the reappearance of the twisted Dirac matrices $\gamma_\mu, \gamma_\mu^*$ of \eqref{eq:dirac} in the above expressions. 

\begin{prop}
\label{prop:instantons}
The instanton gauge potential $\omega=\Psi^* \dd \Psi$ transforms under the action of the Hopf algebra $\U_\theta(so(5,1))$ as $\omega \mapsto \omega+\delta \omega_i$, where
\begin{align*}
&\delta \omega_0 := H_0(\omega)=-z_0 \omega - \half d z_0 \I_2 + \Psi^* ~\gamma_0 ~\dd\Psi,\\
&\delta \omega_1 :=G_{+1,0}(\omega)= -z_1 \omega - \half d z_1 ~\I_2 + \Psi^* ~\gamma_1~ \dd\Psi,\\
&\delta \omega_2 :=G_{0,+1}(\omega)= -z_2 \omega - \half d z_2 ~ \I_2 + \Psi^* ~\gamma_2~ \dd\Psi,\\
&\delta \omega_3 :=G_{-1,0}(\omega)= - \omega \bar z_1-\half d \bar z_1 ~\I_2+\Psi^*~\gamma_1^* ~\dd\Psi,\\
&\delta \omega_4 :=G_{0,-1}(\omega)= -\omega \bar z_2- \half d \bar z_2 ~ \I_2+\Psi^*~\gamma_2^*~ \dd\Psi,
\end{align*}
with $\gamma_\mu, \gamma_\mu^*$ the twisted $4\times 4$ Dirac matrices defined in \eqref{eq:dirac}.
\end{prop}
\begin{proof}
The action of $H_0$ on the gauge potential $\omega=\Psi^* \dd \Psi$ takes the form
\begin{align*}
H_0(\omega) &=H_0(\Psi^*) \dd \Psi + \Psi^* \dd(H_0(\Psi))= \Psi^* 
(- z_0 \I_4 + \gamma_0 )\dd \Psi - \half d z_0 \Psi^* \Psi,
\end{align*}
since $z_0$ is central. Direct computation results in the above expression for $\delta \omega_0$. On the other hand, the twisted action of $G_r$ on 
$\omega=(\omega_{ij})$ takes the form,
\begin{align*}
G_r: \omega_{ij} 
\mapsto \sum_{a,b,c} \tilde\Gamma_{ab} \Psi^*_{ib} (\lambda^{-r_1 H_2})_{ac} \dd \Psi_{cj} +(\lambda^{r_2 H_1})_{ab} \Psi^*_{ib} \Gamma_{ac} \dd \Psi_{cj}  \\ 
\qquad + (\lambda^{r_2 H_1})_{ab} \Psi^*_{ib} (\dd \Gamma_{ac}) \Psi_{cj} ,
\end{align*}
where we used the fact that $\tilde H_j= \sigma H_j \sigma^{-1} = - H_j$.
Let us consider the case $r=(+1,0)$. Firstly, note that the complex numbers $(\lambda^{-H_2})_{ac}$ commute with $\Psi^*_{ib}$ so that from the definition of $\Gamma$ and $\tilde\Gamma$, we obtain for the first two terms,
\begin{align*}
-z_1 (\Psi^* \dd \Psi)_{ij} + 
\half \Psi^*_{ib} 
(\sigma\gamma_1\sigma^{-1})_{cb} (\lambda^{-H_2})_{cd} \dd \Psi_{dj} +\half \Psi^*_{ib} (\gamma_1)_{bc} \dd \Psi_{cj}.
\end{align*}
The first term forms the matrix $-z_1 \omega$ whereas the second two terms combine to give $\Psi^* \gamma_1 \dd \Psi$. 
Finally,  the term $ \Psi^*_{ib} (\dd \Gamma_{ac}) \Psi_{cj}$ reduces to $-\half\dd z_1 \Psi^*_{ib} \Psi_{bj}=-\half \dd z_1 \I_2$. The formulae for $r=(-1,0)$ and $r=(0, \pm 1)$ are established in a similar manner.
\end{proof}

The transformations in Proposition~\ref{prop:instantons} of the gauge potential $\omega$ under the twisted symmetry $\U_\theta(so(5,1))$ induce natural transformations of the canonical connection $\nabla_0$ in \eqref{cancon} 
to $\nabla_{t,i}:=\nabla_0 + t \delta \omega_i + \cO(t^2)$. We shall presently see explicitly that these new connections are (infinitesimal) instantons, {i.e.} their curvatures are self-dual. 
In fact, this also follows from Lemma~\ref{lma:so51-hodge} which states that $\U_\theta(so(5,1))$ acts by conformal transformation therefore leaving invariant the self-duality equations $\ast_\theta F_0 = F_0$ for the basic instanton $\nabla_0$.

We start by writing $\nabla_{t,i}$ in terms of the canonical connection on $\E\isom p \big(\A(\S^4)\big)^4$. Using the explicit isomorphism, giving $\ce$ as equivariant maps $\E\isom \Cinf(\Sk^7) \boxtimes_\rho \C^2$, we find that $\nabla_{t,i}=p \dd + t \delta \alpha_i + \cO(t^2)$ with explicit expressions 
\[
\begin{aligned}
\delta \alpha_0 = p \gamma_0 (\dd p) p - \half \Psi \dd z_0 \Psi^*, \\
\delta \alpha_1 = p \gamma_1 (\dd p) p - \half \Psi \dd z_1 \Psi^*,\\
\delta \alpha_2 = p \gamma_2 (\dd p) p - \half \Psi \dd z_2 \Psi^*,
\end{aligned}
\qquad
\begin{aligned}
~ \\
\delta \alpha_3 = p \gamma_1^* (\dd p) p - \half \Psi \dd z_1^* \Psi^*,\\
\delta \alpha_4 = p \gamma_2^* (\dd p) p - \half \Psi \dd z_2^* \Psi^*,
\end{aligned}
\]
The  $\delta \alpha_i$'s are $4 \times 4$ matrices with entries in the one-forms $\Omega^1(\S^4)$ and  satisfying $p \delta \alpha_i = \delta \alpha_i p = p \delta \alpha_i p = \delta \alpha_i$, as expected from the general theory on connections on modules in 
Appendix~\ref{se:connections}. Indeed, one could move the $\dd z$'s to the left of $\Psi$ at the cost of some $\mu$'s, so getting expression like $(\dd z_i)p \in M_4(\Omega^1(\S^4))$.

From equation~\eqref{ucurv}, the curvature $F_{t,i}$ of the connection $\nabla_{t,i}$ is given by 
\[
F_{t,i}=F_0 + t p \dd (\delta \alpha_i) + \cO(t^2).
\]
To check self-duality  (modulo $t^2$) of this curvature, we will express it in terms of the projection $p$ and consider $F_{t,i}$ as a two-form valued endomorphism on $\E$. 
\begin{prop}
\label{prop:curv-transf}
The curvatures $F_{t,i}$ of the connections $\nabla_{t,i}$, $i=0,\ldots,4$, are given by $F_{t,i} =F_0 + t \delta F_i+ \cO(t^2)$, where $F_0=p \dd p \dd p$ and 
\[
\begin{aligned}
\delta F_0&=-2 z_0 F_0, \\
\delta F_1&=-2 z_1 \lambda^{H_2} F_0, \\
\delta F_2&=-2 z_2 \lambda^{H_1} F_0;
\end{aligned}
\qquad
\begin{aligned}
~ \\
\delta F_3=-2 z_1^*\lambda^{-H_2} F_0, \\
\delta F_4=-2 z_2^*\lambda^{-H_1} F_0.
\end{aligned}
\]
\end{prop}
\begin{proof}
A small computation yields for $\delta F_i=p \dd (\delta \alpha_i)$, as an $\Omega^2(\S^4)$-valued  endomorphism on $\E$ the expression, 
$\delta F_i = p (\dd p) \gamma_i (\dd p) p - p \gamma_i (\dd p) (\dd p) p$,
with the notation $\gamma_3=\gamma_1^*$ and $\gamma_4=\gamma_2^*$, and using $p (\dd p) p=0$. 
Then, the crucial property $p (\dd p\gamma_i +\gamma_i \dd p) (\dd p) p = 0$ all $i=0,\ldots,4$ yields  $\delta F_i = -2 p \gamma_i \dd p \dd p p$. This is expressed as $\delta F_i = -2 p \gamma_i p \dd p \dd p$ by using 
$\dd p = (\dd p) p + p \dd p $. Finally, $p \gamma_i p=\Psi (\Psi^* \gamma_i \Psi) \Psi^*$, and the result follows from the definition of the $z$'s in terms of the Dirac matrices given in \eqref{eq:dirac}, together with their commutation relations with the matrix $\Psi$.
\end{proof}
\begin{prop}
The connections $\nabla_{t,i}$ are (infinitesimal) instantons, {i.e.} 
\[
\ast_\theta F_{t,i}=F_{t,i}  \qquad \mathrm{mod} \, t^2 .
\]
Moreover, the
connections $\nabla_{t,i}$ are not gauge equivalent to $\nabla_0$.\end{prop}
\begin{proof}
The first point follows directly from the above expressions for $\delta F_i$ and the self-duality of $F_0$.   
To establish the not gauge equivalence, recall from \eqref{infgt} that an infinitesimal gauge transformation is $\nabla_0 \mapsto \nabla_0 + t [\nabla_0,X]$ for $X \in \Gamma(\ad(\Sk^7))$. We need to show that $\delta \omega_i$ is orthogonal to $[\nabla_0,X]$ for any such $X$, {i.e.}
that 
$$
([\nabla_0,X],\delta \omega_i )_2= 0,
$$ 
with the natural inner product on $\Omega^1(\ad(\Sk^7):=\Omega^1(\S^4)\ot_{\Cinf(\S^4)} \Gamma(\ad(\Sk^7))$. Now, 
one has that 
$$(\nabla_0^{(2)}(X),\delta \omega_i)_2 = (X, \left(\nabla_0^{(2)}\right)^* (\delta \omega_i) )_2,$$ 
which then should vanish for all $X$. But  $\delta \omega_i = T_i(\omega)$ coincides with $L_\theta (t_i \cdot \omega^\class)$ with $t_i$ and $\omega^\class$ the classical counterparts of $T_i$ and $\omega$, respectively. In the undeformed case, the infinitesimal gauge potentials generated by acting with elements in $so(5,1)-so(5)$ on the basic instanton gauge potential $\omega^\class$ satisfy $(\nabla_0^{(2)} )^* (\delta \omega^\class_i)=0$ as shown in \cite{AHS78}. The result then follows from the observation that $\nabla_0^{(2)}$ commutes with the quantization map $L_\theta$.
\end{proof}

In \cite{LS06} we gave a 
completeness argument on the family of instantons presented above, by index theoretical arguments, similar to the one in \cite{AHS78} for undeformed instantons on $S^4$. The dimension of the `tangent' of the moduli space can be computed as the index of a twisted Dirac operator which turns out to be equal to its classical value that is five.

\section{An instanton bundle over a symplectic quantum sphere $S^4_q$}

A different quantum version of  the $SU(2)$ Hopf bundle $S^7 \rightarrow S^4$ was constructed in \cite{LPR06}. The quantum sphere $S^7_q$  arises from the symplectic group $Sp_q(2)$ and a quantum $4$-sphere $S^4_q$ is obtained via a suitable self-adjoint idempotent $p$ whose entries generate the algebra
$A(S^4_q)$ of polynomial functions over it -- a property in common with the toric four sphere $\S^4$ used above. This projection determines a
deformation of an instanton bundle over the classical sphere $S^4$.

One starts with the symplectic quantum groups $A(Sp_q(2))$, 
i.e. the Hopf algebras generated by matrix elements
$T_i^j$'s with commutation rules coming  from the
$R$ matrix of the $C$-series \cite{frt}. The symplectic quantum
$7$-sphere $A(S_q^7)$ is generated by the matrix elements of
the first and
the last  columns of $T$. The
algebra $A(S_q^7)$ is the  quantum version of the homogeneous space
$Sp(2)/Sp(1)$
and the injection
$A(S^7_q)\hookrightarrow A(Sp_q(2))$ is a quantum principal bundle with `structure Hopf algebra'  $A(Sp_q(1))$,
an example of the general construction of \cite{BM93}.

In turn, the sphere $S_q^7$ is the total space of a
quantum $\SU_q(2)$
principal bundle over a quantum $4$-sphere $S_q^4$. 
Unlike the construction for $A(S_q^7)$ out of $A(Sp_q(2))$, now one does not have a quantum homogeneous structure. Still, there is a natural coaction of
$\SU_q(2)$ on $A(S^7_q)$ with coinvariant algebra $A(S^4_q)$ and the injection $A(S^4_q)\hookrightarrow A(S^7_q)$ is another instance of a quantum principal bundle. 

To compute the charge of the projection $p$ and to prove the non triviality of this principal bundle, one follows a general strategy of noncommutative index theorem \cite{C94}. One constructs the 
representations of the algebra $A(S_q^4)$ and the corresponding 
$K$-homology. 
The
analogue of the fundamental class of $S^4$ is given by a non trivial Fredholm
module $\mu$. The natural coupling between $\mu$ and the projection 
$p$ is computed
via the pairing of the corresponding Chern characters $\chern ^*(\mu)\in
HC^*[A(S_q^4)]$ and
$\chern_*(p)\in HC_*[A(S_q^4)]$ in cyclic cohomology and homology respectively. 
The result of this pairing, which is an integer by
principle -- being the index of a Fredholm operator -- is indeed $-1$ and the bundle is non trivial.

Clearly the next step would be to repeat the analysis of the toric four sphere and define a Yang--Mills action functional  and self-duality equations. To this end one needs a `metric structure' on the bundle; for this, the recently found \cite{DLSSV04} isospectral noncommutative geometry for  $\SU_q(2)$ promises to be useful.

\subsection{Odd spheres from quantum symplectic groups}
There are quantum spheres associated with the 
compact real form
of the quantum symplectic groups
$Sp_q(N, \IC)\; (N=2n)$, the latter being given in
\cite{frt}. 
The algebra $A(Sp_q(N, \IC))$ is the associative noncommutative 
algebra generated
over the ring of Laurent polynomials $\IC_q:=\IC [q,q^{-1}]$ by the entries
${T_i}^j,\; i,j=1,\dots,N$ of a  matrix $T$ which satisfy RTT equations:
$$ 
R\;T_1 T_2 = T_2 T_1 R~, \qquad  \quad T_1 := T \ot 1 \; ,
\quad T_2 := 1 \ot T \; .
$$ 
Here the
relevant $N^2 \times N^2$ matrix $R$ is the one for the  $C_N$ series 
and has the
form
\cite{frt},
\begin{align*} 
R &= q\sum_{i=1}^N {e_i}^i
\ot {e_i}^i +
\sum_{\stackrel{i,j=1}{i \neq j,{j'}}}^N {e_i}^i \ot {e_j}^j + q^{-1}
\sum_{i=1}^N {e_{i'}}^{i'} \ot {e_i}^i
\\ & \qquad + (q -q^{-1}) \sum_{\stackrel{i,j=1}{i >j}}^N {e_i}^j
\ot {e_j}^i -(q -q^{-1})
\sum_{\stackrel{i,j=1}{i >j}}^N q^{\rho_i - \rho_j} \varepsilon_i
\varepsilon_j ~{e_i}^j \ot {e_{i'}}^{j'}  ~,
\end{align*} 
where:
\begin{itemize}
\item[ ] $i'= N+1 -i$ ~;
\item[ ] ${e_i}^j \in M_n(\IC)$  are the elementary  matrices, i.e.
$({e_j}^i)^k_l=\delta_{jl}\delta^{ik}$ ~;
\item[ ] $\varepsilon_i =1, \mbox{ for } i=1,\dots, n$ ~;
\quad  $\varepsilon_i = -1 , \mbox{ for }i=n+1 ,\dots, N$ ~;
\item[ ] $(\rho_1,\dots,\rho_N) =(n,n-1,\dots , 1,-1, \dots ,-n)$ ~.
\end{itemize}

\bigskip
The symplectic group structure comes from the matrix
${C_i}^j=q^{\rho_j}\varepsilon_i \delta_{ij'}$ by  imposing  the additional
relations
$$ 
TCT^tC^{-1}=CT^tC^{-1}T=1 \;.
$$ 
The Hopf algebra structure $(\Delta, \varepsilon, S)$ of the 
quantum group
$Sp_q(N,\IC)$ is given by
$$
\Delta (T) =T \pot T \; ,  \quad
\varepsilon(T)=I \; , \quad  S(T)=CT^tC^{-1} \; .
$$ In components the antipode explicitly reads
\[
{S(T)_i}^j = -q^{\rho_{i'}+\rho_j} \varepsilon_i
\varepsilon_{j'}{T_{j'}}^{i'} \; .
\] 
At $q=1$ the Hopf algebra  $Sp_q(N,\IC)$ reduces to 
the algebra of
polynomial functions over the symplectic group  $Sp(N,\IC)$.
The compact real form $A(Sp_q(n))$ is given by 
taking $q\in \IR$
and the natural anti-involution
\begin{equation}\label{inv}
\overline{T}=S(T)^t=C^t T (C^{-1})^t \; .
\end{equation}

Let us denote
$$ 
x_i={T_i}^N \;, \quad v^j = {S(T)_N}^j \; , \quad i,j = 1, \dots, N \; .
$$ 
The elements $\{x_i,~v^j\}$ close an algebra as we show presently.
By summing over repeated indexes, the  RTT equations in components read 
\begin{equation}\label{rtt} 
{R_{ij}}^{kp} ~{T_k}^r ~{T_p}^s = {T_j}^p ~{T_i}^m
~{R_{mp}}^{rs} \; .
\end{equation} 
From which, some algebra yields the relations,
\be\label{st}
\begin{aligned}
& {S(T)_a}^i~{S(T)_p}^j~{R_{ij}}^{kl} = {R_{ap}}^{rs}
~{S(T)_s}^l~{S(T)_r}^k \;,  \\
& {S(T)_l}^j ~{R_{ij}}^{kp}~{T_k}^r ={T_i}^m ~{R_{ml}}^{rs} ~{S(T)_s}^p \; .  
\end{aligned}
\ee

From \eqref{rtt} with $r=s=N$ we have
\[ 
{R_{ij}}^{kp}~ x_k x_p = {T_j}^p ~{T_i}^m ~{R_{mp}}^{NN} \;.
\] 
Since the only element ${R_{mp}}^{NN} ~ \varpropto ~{e_m}^N \ot
{e_p}^N  ~(m,p\leq N)$  which is different from zero is ${R_{NN}}^{NN}=q$,  it
follows that
\begin{equation}\label{comxx}
    {R_{ij}}^{kp}~ x_k x_p = q ~x_j x_i \; .
\end{equation} 

Putting  $a=p=N$ in the first of equations \eqref{st}, we get
$$ v^i v^j {R_{ij}}^{kl}= {R_{NN}}^{rs} {S(T)_s}^l {S(T)_r}^k \; .
$$ The sum on the r.h.s. reduces to ${R_{NN}}^{NN} {S(T)_N}^l 
{S(T)_N}^k$ and the
$v^i$'s give an algebra with commutation relations
\begin{equation}\label{commvv} v^l v^k {R_{lk}}^{ji} = q v^i v^j .
\end{equation} 

Finally, for  $l=r=N$ the second of equations \eqref{st} reads:
    $$ v^j {R_{ij}}^{kp} ~x_k = {T_i}^m {R_{mN}}^{Ns} {S(T)_s}^p.$$ 
The
only term in $R$ of the form ${e_m}^N \ot {e_N}^s ~(m\leq N)$  is  ${e_N}^N \ot
{e_N}^N$; thus 
\begin{equation}\label{commxv} v^j {R_{ij}}^{kp} ~x_k = q ~x_i v^p  \; .
\end{equation} 

\bigskip
With the anti-involution \eqref{inv} we have 
$v^j={S(T)_N}^i= \bar{x}^j$, 
and the subalgebra of
$A(Sp_q(n))$ generated by
$\{x_i, v^i=\bar{x}^i, ~i=1,\dots,2n\}$ is the algebra $A(S^{4n-1}_q)$ of polynomial 
functions on a
sphere. Indeed
$$  
S(T)T=I~ \Rightarrow  ~\sum {S(T)_N}^i {T_i}^N =\delta_N^N =1 \;,
$$  
from which we get a sphere relation
\[
\sum_i \bar{x}^i x_i = 1 \; .
\] 
Furthermore, the restriction of the comultiplication 
is a natural
left coaction
$$
\Delta_L : A(S^{4n-1}_q) \lr A(Sp_q(n)) \ot A(S^{4n-1}_q)\; .
$$ The fact that $\Delta_L$ is an algebra map then implies that 
$A(S^{4n-1}_q)$  is
a comodule algebra over
$A(Sp_q(n))$.

At $q=1$ this algebra is made of polynomial functions over the
spheres $S^{4n-1}$ as homogeneous  spaces of the symplectic group $Sp(n):\;
S^{4n-1}=Sp(n)/Sp(n-1)$.

\subsection{The symplectic quantum sphere $S^7_q$}\label{s7}

The algebra
$A(S^7_q)$ is generated by the elements $x_i = {T_i}^4$ and
$\bar{x}^i={S(T)_4}^i= q^{2+\rho_i}
\varepsilon_{i'}{T_{i'}}^1,$ for $i=1,\dots ,4$ with sphere relation $\sum_{i=1}^4 \bar{x}^i x_i = 1$
coming from $S(T)~T = 1$. We explicitly give the commutation relations among generators  systematically used in the following. \\
   From \eqref{comxx}, the algebra of the $x_i$'s is given by
\[
\begin{array}{ll}
\uno \due = q \due \uno \; , & \uno \tre = q \tre \uno \; ,
\\
\due \qu = q \qu \due \; , & \tre \qu = q \qu \tre \; ,
\\
\qu \uno = q^{-2} \uno \qu \; , &
\tre \due = q^{-2} \due \tre +q^{-2}(q^{-1}-q) \uno \qu \; ,
\end{array}
\] 
together with their conjugates (given in \eqref{commvv}). The commutation relations between the
$x_i$ and the $\bar{x}^j$ are deduced  from \eqref{commxv}:
\[
\begin{array}{l}
\begin{array}{ll}
\uno  \unob = \unob \uno \; ,&
\uno \dueb = q^{-1} \dueb \uno \; ,
\\
\uno \treb = q^{-1} \treb \uno \; , &
\uno \qub = q^{-2} \qub \uno \; ,
\end{array}
\\
\\
\begin{array}{l}
\due \dueb = \dueb \due + (1-q^{-2})\unob \uno \; ,\\
\due \treb = q^{-2} \treb \due \; , \\
\due \qub = q^{-1} \qub \due +q^{-1} (q^{-2}-1) \treb \uno \; ,
\end{array}
\\
\\
\begin{array}{ll}
\tre \treb = \treb \tre +(1-q^{-2})[\unob \uno +(1+q^{-2})\dueb \due] \; ,
\\
\tre \qub = q^{-1} \qub \tre +(1-q^{-2})q^{-3} \dueb \uno \; ,
\end{array}
\\
\\
\begin{array}{ll}
\qu \qub = \qub \qu +(1-q^{-2})[(1+q^{-4})\unob \uno + \dueb \due + \treb
\tre] \; ,
\end{array}
\end{array}
\] 
again with their conjugates.

In complete analogy with  the classical homogeneous space
$Sp(2)/Sp(1) \simeq S^7$,
the algebra $A(S^7_q)$ can be realized as the subalgebra of
$A(Sp_q(2))$ made of elements which are coinvariant  under the right-coaction of
$A(Sp_q(1))$.

\begin{lma} The two-sided *-ideal in $A(Sp_q(2))$ generated by
$$ 
I_q = \{{T_1}^1-1, {T_4}^4-1, {T_1}^2, {T_1}^3, {T_1}^4, {T_2}^1, {T_2}^4, {T_3}^1, {T_3}^4, {T_4}^1, {T_4}^2,{T_4}^3\} \; ,
$$ 
with the involution \eqref{inv} is a Hopf ideal.
\end{lma}
\noindent {\it Proof.} Since  ${S(T)_i}^j \propto{T_{j'}}^{i'}$,
$S(I_q)\subseteq
I_q$ which also proves that $I_q$ is a *-ideal.  One easily  shows that
$\varepsilon (I_q)=0$ and $\Delta (I_q)\subseteq I_q\ot
A(Sp_q(2))+A(Sp_q(2))\ot
I_q$. \qed

\begin{prop} The Hopf algebra $B_q := A(Sp_q(2))/I_q$  is isomorphic to the
coordinate algebra $A(SU_{q^2}(2)) \cong A(Sp_q(1))$.
\end{prop}
\noindent {\it Proof.} Using $\overline{T}=S(T)^t$ and setting
${T_2}^2 = \alpha,~
{T_3}^2=\gamma $, the algebra $B_q$ can be described as the algebra
generated by the
entries of the matrix
\[ 
T'=\left(
\begin{array}{cccc} 1& 0 & 0 &0
\\ 0 & \alpha & -q^2 \bar{\gamma}  & 0
\\ 0 & \gamma &  \bar{\alpha} & 0
\\ 0& 0 & 0 & 1
\end{array}
\right).
\] 
The commutation relations deduced from RTT equations \eqref{rtt}
read:
\[
\begin{array}{ll}
\alpha \bar\gamma = q^2 \bar\gamma \alpha \quad , \quad &
\alpha \gamma = q^2 \gamma \alpha ~,
\quad
\gamma \bar\gamma= \bar\gamma \gamma ~,
\\
\bar\alpha \alpha + \bar\gamma \gamma =1  \quad ;&
\alpha \bar\alpha + q^4 \gamma \bar\gamma =1 ~.
\end{array}
\] 
Hence, as an algebra $B_q$  is isomorphic to the algebra
$A(SU_{q^2}(2))$.  Furthermore, the restriction of the coproduct of 
$A(Sp_q(2))$
to $B_q$  endows the latter with a  coalgebra structure, $\Delta(T') 
= T' \pot T'$,
which is the same as the one of
$A(SU_{q^2}(2))$.  We can conclude that also as a Hopf algebra $B_q$  is
isomorphic to the  Hopf algebra
$A(SU_{q^2}(2)) \cong A(Sp_q(1))$.
\qed

\begin{prop}\label{p2} The algebra  $A(S_q^7) \subset  A(Sp_q(2))$ is
the algebra
of coinvariants with respect to the natural right coaction
$$
\Delta_R:A(Sp_q(2)) \to A(Sp_q(2)) \pot
A(Sp_{q}(1))
\quad , \quad
\Delta_R(T) = T \pot T' \; .
$$
\end{prop}
\noindent
{\it Proof.} It is straightforward to show that the generators of the
algebra $A(S_q^7)$ are  coinvariants:
$$
\Delta_R (x_i)= \Delta_R (T_i^4)= x_i \ot 1 ~~;~~
\Delta_R (\bar{x}^i)= -q^{2 + \rho_i} \varepsilon_i ~\Delta_R (T_i^1)=
\bar{x}^i \ot 1 \,
$$
thus the algebra $A(S_q^7)$ is made of coinvariants. There are no other coinvariants of degree one
since each row of the submatrix of $T$ made out of the two central columns is a fundamental comodule under the coaction of $SU_{q^2}(2)$.
Other coinvariants arising at higher even degree are of the form 
$(T_{i2}T_{i3}-q^2 T_{i3}T_{i2})^n$; the commutation relations of $A(Sp_q(2))$
yield that these belong to $A(S^7_q)$ as well.
On the other hand, similar expressions involving elements from different rows cannot be coinvariant.
\qed \\

The previous construction is one more example of the general construction of a quantum principal bundle over a quantum homogeneous space \cite{BM93}. 
The latter is the datum of a Hopf quotient $\pi: A(G) \to A(K)$ with the right coaction of $A(K)$ on $A(G)$ given by the reduced coproduct 
$\Delta_R:=(id \ot \pi)\Delta$ where $\Delta$ is the coproduct of $A(G)$. The subalgebra $B \subset A(G)$ made of the coinvariant elements with respect
to $\Delta_R$ is called a quantum homogeneous space. To prove that it is a quantum principal bundle one needs some more assumptions (see
Lemma 5.2 of \cite{BM93}). In our case $A(G)=A(Sp_q(2)),~A(K)=A(Sp_q(1))$ with $\pi(T)=T'$. The resulting inclusion $ B=A(S^7_q) \hookrightarrow A(Sp_q(2))$ is indeed a quantum principal bundle.

\subsection{The quantum sphere $S^{4}_q$}

The next  step is to make the sphere $S_q^7$
itself into the total space of a quantum principal bundle over a
deformed $4$-sphere. Unlike what we saw in the previous section, this is not a quantum homogeneous
space construction and it is not obvious that such a bundle exists at all.
Nonetheless the notion of quantum bundle is more general and one only needs
that the total space algebra is a comodule algebra over a Hopf algebra with additional suitable properties.
The precise sense in which the algebra inclusion $A(S^{4}_q) \hookrightarrow  A(S^{7}_q)$ to be constructed is a noncommutative principal bundle was explained in \cite{LPR06} to which we refer for more details. Presently we shall illustrate the main ingredients.

\bigskip
On the free module $\ce:=\IC^4 \ot A(S_q^7)$ we have the natural 
Hermitean structure $\hs{\xi_1}{\xi_2} =  \sum_{j=1}^4 \bar{\xi_1}^j \xi_2^j$, and
to every $\ket{\xi}\in\ce$ one associates an element 
$\bra{\xi}$ in the dual module $\ce^*$ by the pairing,
$$
\bra{\xi}(\ket{\eta}):=\hs{\xi}{\eta}=h(\ket{\xi},\ket{\eta})\;, \quad \forall ~ \eta\in\ce \; .
$$ 
Guided by
the classical construction \cite{Lnd00}, we
look for two elements $\ket{\phi_1},~\ket{\phi_2}$  in $\ce$ with the
property that
$$
\hs{\phi_1}{\phi_1} =1 \;, \quad \hs{\phi_2}{\phi_2} =1 \;, \quad
\hs{\phi_1}{\phi_2} =0 \; .
$$ 
Then, one gets a self-adjoint idempotent (a projection) by the matrix-valued function
$p := \ket{\phi_1} \bra{\phi_1} + \ket{\phi_2}\bra{\phi_2}$.
In 
principle, $ p
\in \Mat_4(A(S^7_q))$, but we choose the pair
$(\ket{\phi_1},~\ket{\phi_2})$ in such a way that the entries of $p$ 
will generate a
subalgebra $A(S^4_q)$ of
$A(S^7_q)$ which is a deformation of the algebra of polynomial functions on the
$4$-sphere $S^4$. The appropriate $\ket{\phi_1},~\ket{\phi_2}$ turn out to be the following:
\bean
&& \ket{\phi_1}= (q^{-3} \uno ,- q^{-1} \dueb , q^{-1} \tre ,-\qub )^t \; ,   \\
&& \ket{\phi_2}= ( q^{-2} \due ,  q^{-1} \unob ,  - \qu ,- \treb )^t \; .  
\eean 
The projection $p$ can be equivalently written as
\begin{equation}\label{u}
p=\Psi \Psi^* \;, \quad \mathrm{with} \quad \Psi=\left( \ket{\phi_1},  \ket{\phi_2}
\right) = \left(
\begin{array}{cc} q^{-3} \uno & q^{-2} \due
\\ -q^{-1}\dueb & q^{-1} \unob
\\
    q^{-1} \tre & -\qu
\\
    -\qub & -\treb
\end{array}
\right) \; .
\end{equation} 
having $\Psi^*\Psi=1$, as can be readily established. 

\begin{prop} The entries of the projection
$p=\Psi \Psi^*$ in \eqref{u}, generate a subalgebra of $A(S^7_q)$ 
which is a
deformation of the algebra of polynomial functions on the
$4$-sphere $S^4$. Esplicitly, 
\begin{equation}\label{proj} p=\left(
\begin{array}{llll} q^{-2}t & 0 & a &b
\\
\\ 0 &t & q^{-2} \bar{b} & -q^2 \bar{a}
\\
\\
\bar{a} & q^{-2} b & 1-q^{-4}t & 0
\\
\\
\bar{b}& -q^2 a & 0 & 1-q^2t
\end{array}
\right) \; ,
\end{equation} 
and generators given by
\bean 
&& t = q^{-2} \dueb \due  + q^{-2} \unob \uno \; , \nn \\
&& a=q^{-4} \uno \treb - q^{-2} \due \qub \; , \qquad 
b =-q^{-3} \uno \qu -q^{-2} \due \tre
\eean
\end{prop}
\noindent {\it Proof.} The above fourmul{\ae} are obtained by direct computation. 
By construction $p^*=p$ and this means that 
$\bar{t}=t$, 
and that $\bar{a},\bar{b}$ are conjugate to $a,b$ respectively.   Also, by
construction $p^2=p$ and  this gives the easiest way to compute the 
commutation relations
between the generators. One finds sphere relations
\begin{equation}\label{sr4}
\begin{array}{lll} a\bar{a}+b\bar{b}=q^{-2}t(1-q^{-2}t) \; , ~~~~~ & 
q^4 \bar{a}a
+q^{-4}\bar{b}b=t(1-t) \; ,
\\ b \bar b - q^{-4}\bar{b}b =(1- q^{-4})t^2\; .&
\end{array}
\end{equation}
and 
\begin{equation}\label{s4}
\begin{array}{ll} ab = q^4 ba \; , & \; \bar{a}b =b \bar{a} \; , \\ 
ta= q^{-2} at
\; , & \;tb = q^4 bt \; ,
\end{array}
\end{equation}
together with their conjugates. \qed

 We define the algebra $A(S^4_q)$ to be  the algebra generated by the
elements $a,\bar{a},b,\bar{b},t$ with the commutation relations \eqref{s4} and
\eqref{sr4}. At $q=1$ it reduces to the algebra of polynomial functions on the
sphere $S^4$. We can take $|q|<1$, since
$$q\mapsto q^{-1},\quad a\mapsto q^2\bar a,\quad b\mapsto q^{-2}\bar b,
\quad t\mapsto q^{-2} t$$ yields an isomorphic algebra. \\

\subsection{The principal bundle  $A(S^{4}_q) \hookrightarrow  A(S^{7}_q)$}
There is a coaction of the quantum group $SU_q(2)$ on 
the sphere
$S^7_q$ that makes up  the quantum principle bundle structure.
Notice that the two pairs of generators $(\uno,\due),
(\tre,\qu)$ both yield a quantum plane:
\begin{eqnarray*}
\uno \due = q \due \uno \; ,  & \quad& \unob \dueb = q^{-1} \dueb \unob \; ,
\\
\tre \qu = q \qu \tre \; , && \treb \qub =q^{-1} \qub \treb \; .
\end{eqnarray*} Then we shall look for a right-coaction of $SU_q(2)$ 
on the rows of
the matrix $\Psi$ in \eqref{u}. Other pairs of generators yield
quantum planes but
the only choice which gives a projection with the right number of 
generators is the one given above.

The  defining matrix of the quantum group $SU_q(2)$ reads \cite{W87}
\[
\begin{pmatrix}
\alpha  & -q \bar{\gamma} \\
\gamma & \bar{\alpha}
\end{pmatrix}
\]
with commutation relations,
\begin{equation}\label{SU2}
\begin{array}{lll}
\alpha \gamma = q \gamma \alpha \; , \qquad & \alpha \bar{\gamma} = q
\bar{\gamma} \alpha \; , \qquad & \gamma 
\bar{\gamma}=\bar{\gamma}\gamma \; ,  \\
\alpha \bar{\alpha} + q^2 \bar{\gamma} \gamma = 1 \; , & \bar{\alpha} \alpha +
\bar{\gamma} \gamma = 1 \; .
\end{array}
\end{equation} We define a coaction of $SU_q(2)$ on the matrix \eqref{u} by,
\begin{equation}\label{cosu2}
\delta_R (\Psi) :=
\left(
\begin{array}{cc} q^{-3} \uno & q^{-2} \due
\\- q^{-1}\dueb & q^{-1} \unob
\\
    q^{-1} \tre & -\qu
\\ - \qub & -\treb
\end{array}
\right)
\stackrel{.}{\otimes}
\begin{pmatrix}
\alpha  & -q \bar{\gamma} \\
\gamma & \bar{\alpha}
\end{pmatrix} \; .
\end{equation} We shall prove presently that this coaction comes from 
a coaction of
$A(SU_q(2))$ on the  algebra $A(S^7_q)$. For the moment we 
remark that, by
its form in
\eqref{cosu2} the entries of the projection $p=\Psi \Psi^*$ are automatically
coinvariants.

On the generators, the coaction \eqref{cosu2} is given explicitly by
\begin{equation} \label{coactionp}
\begin{array}{ll}
\delta_R(\uno)=  \uno \otimes \alpha + q ~\due \otimes \gamma \; ,  &
\delta_R(\unob) = q \dueb \otimes \bar{\gamma} + \unob \otimes
\bar{\alpha} =\overline{\delta_R(\uno)} \; ,
\\
\delta_R(\due)= - \uno \otimes \bar{\gamma} + \due \otimes \bar{\alpha} \; ,  &
\delta_R(\dueb) = \dueb \otimes \alpha - \unob \otimes \gamma
=\overline{\delta_R(\due)} \; ,
\\
\delta_R(\tre)=  \tre \otimes \alpha - q~\qu \otimes \gamma \; ,  &
\delta_R(\treb) =- q \qub \otimes \bar{\gamma} + \treb \otimes
\bar{\alpha}  =\overline{\delta_R(\tre)} \; ,
\\
\delta_R(\qu)=  \tre \otimes \bar{\gamma} + \qu \otimes \bar{\alpha} \; ,  &
\delta_R(\qub) = \qub \otimes \alpha + \treb \otimes  \gamma 
=\overline{\delta_R(\qu)}
\; ,
\end{array}
\end{equation} 
from which it is clear its compatibility with the
anti-involution, i.e. $\delta_R(\bar{x}^i)=\overline{\delta_R(x_i)}$. The 
map $\delta_R$  extends as an algebra homomorphism to the whole of
$A(S^{7}_q)$. We have the following:
\begin{prop}\label{p4} The coaction \eqref{coactionp} is a right coaction of the quantum
group $SU_q(2)$
on the 7-sphere  $S^7_q$,
\[
\delta_R : A(S^{7}_q) \lr A(S^{7}_q) \ot A(SU_q(2)) \; .
\]
Moreover, the algebra $A(S_q^4)$ is the algebra of corresponding
coinvariant elements.

\end{prop}
\noindent
\textit{Proof}.
For the first point, by using the commutation relations of $A(SU_q(2))$ in
\eqref{SU2}, a lengthy computation gives that the 
commutation relations of
$A(S^7_q)$ are preserved under the map $\delta_R$. This also shows that extending $\delta_R$ as an
algebra homomorphism yields a consistent coaction.

Next, one has to show that 
$A(S^4_q)= \{ f \in A(S^7_q) ~|~ \delta_R(f)=f \ot 1 \}$.
By using the commutation relations of $A(S_q^7)$ and those of
$A(SU_q(2))$, one prove also by direct calculation that the generators of
$A(S_q^4)$ are coinvariants, thus showing that $A(S_q^4)$ is made of coinvariants (we have already observed that the elements of $p$ are automatically coinvariant).
There are no other coinvariants for the following reason.
From equation \eqref{coactionp} it is clear that $w_1 \in \{\uno, \tre, \dueb, \qub \}$ (respectively $w_{-1} \in \{ \due, \qu, \unob, \treb \}$) are weight vectors of weight $1$
(resp. $-1$) in the fundamental comodule of $SU_q(2)$.
It follows that the only possible coinvariants are of the form $(w_1 w_{-1} -q w_{-1} w_1)^n$. When $n=1$ these are just the generators of
$A(S^4_q)$. \qed

\noindent The right coaction of $SU_q(2)$ on the 7-sphere  $S^7_q$ 
can be written as
\[
\delta_R(\uno , \due ,\tre ,\qu)= (\uno , \due , \tre , \qu)
\stackrel{.}{\otimes}
\left(
\begin{array}{cccc}
\alpha ~~& -\bar{\gamma} ~~ &0 ~~ &0
\\ q \gamma & \bar{\alpha} & 0 &0
\\ 0 & 0 & \alpha &  \bar{\gamma}
\\ 0 &0 &-q \gamma &\bar{\alpha}
\end{array}
\right) \; ,
\] 
together with $\delta_R(\bar{x}_i)= 
\overline{\delta_R(x_i)}$. In the
block-diagonal matrix which appears in the equation above, the first 
copy is given
by $SU_q(2)$ while the second one is twisted as
$$
\begin{pmatrix}
\alpha & \bar{\gamma} \\ -q \gamma &  \bar{\alpha}
\end{pmatrix}
    =
\begin{pmatrix} 1 & 0 \\ 0& - 1
\end{pmatrix}
\begin{pmatrix}
\alpha & -\bar{\gamma} \\ q \gamma &  \bar{\alpha}
\end{pmatrix}
\begin{pmatrix} 1 & 0 \\ 0& - 1
\end{pmatrix} \; .
$$
We refer to \cite{LPR06} for more details on the structure of the algebra inclusion $\A(\S^4) \into \A(\Sk^7)$ as a noncommutative principal bundle. 

\subsection{The index pairings}\label{in-pa}

The `defining' self-adjoint idempotent $p$ in \eqref{proj} determines 
a class in
the $K$-theory of $S_q^4$, i.e. $[p]\in K_0[\cc(S_q^4)]$. A way to establish its
nontriviality is by pairing it with a nontrivial element in the dual 
$K$-homology,
that is with (the class of) a nontrivial Fredholm module $[\mu]\in 
K^0[\cc(S_q^4)]$.
In fact,  to compute the pairing, it is more convenient to first compute the
corresponding Chern characters in the cyclic homology
$\chern_*(p) \in HC_*[A(S^{4}_q)]$ and cyclic cohomology $\chern^*(\mu)\in
HC^*[A(S^{4}_q)]$ respectively, and then use the pairing between 
cyclic homology
and cohomology \cite{C94}. To compute the pairing and to prove the nontriviality of the
bundle it is enough to consider $HC_0[A(S^{4}_q)]$ and dually to take a suitable trace of the projector. 

\bigskip
The Chern character of the projection $p$ in \eqref{proj} has a component in degree zero,
$\chern_0(p)\in HC_0[A(S^{4}_q)]$, simply given by the matrix trace,
\[
\chern_0(p) := \tr(p) = 2 - q^{-4} (1-q^2)(1-q^4) ~t ~\in A(S^4_q).
\] 
The higher degree parts of $\chern_*(p)$ are obtained via the
periodicity operator $S$.

\bigskip
The K-homology of an involutive algebra $A$  is given in
terms of homotopy classes of Fredholm modules. For this we need representations of the algebra. 

\subsection{Representations of the algebra $A(S^4_q)$}\label{se:rep} 

We construct irreducible $*$-representations  of $A(S^4_q)$ as bounded operators on a separable Hilbert space $\hil$. Since $q\mapsto q^{-1}$
gives an isomorphic algebra, we restrict ourselves to  $|q|<1$. 

We consider representations which are $t$-finite \cite{KS97}, i.e. such that the
eigenvectors of $t$ span $\hil$. Since the self-adjoint operator
$t$ must be bounded due to the spherical relations, from the 
commutation relations
$ta= q^{-2} at,~~ t\bar{b}=q^{-4}\bar{b}t,$ it follows that the 
spectrum should be
of the form $\lambda q^{2k}$ and
$a, \bar b$ (resp. $\bar a, b$) act as rising (resp. lowering) operators on the
eigenvectors of $t$. Then boundedness implies the existence of a 
highest weight
vector, i.e. there exists a vector
$\ket{0,0}$ such that
\[
t\ket{0,0}=t_{00}\ket{0,0},~~~a\ket{0,0}=0,~~\bar{b}\ket{0,0}=0 ~.
\]
By evaluating $q^4 \bar{a} a + b \bar{b}=(1-q^{-4}t)t$ on $\ket{0,0}$ we have
$$
(1-q^{-4}t_{00})t_{00}=0 \;.
$$ 
According to the values of the eigenvalue $t_{00}$ we have two representations.
 
\subsubsection{The representation $\beta$}\label{se:pi} 
The first 
representation,
that we call $\beta,$ is obtained for
$t_{00}=0$. Then, $t\ket{0,0}=0$ implies
$t=0$. Moreover, using the commutation relations \eqref{s4} and \eqref{sr4}, it
follows that this representation  is the trivial one
\begin{equation} \label{pi} t=0,~~~a=0,~~~b=0 \; ,
\end{equation} the representation Hilbert space being just $\IC$; of 
course, $\beta
(1)=1$.

\subsubsection{The representation $\sigma$}\label{se:sigma} 
The second
representation, that we call $\sigma,$ is obtained for
$t_{00}=q^4$ and  is infinite dimensional. An
orthonormal basis of the representation Hilbert space $\hil$
is given by the set
$\ket{m,n}= N_{mn} \bar{a}^m b^n \ket{0,0}$, with $n,m \in \IN$, 
$N_{00}=1$ and
$N_{mn}\in\IR$ normalization constants. Then
$$
\begin{array}{l} t\ket{m,n}=t_{mn}\ket{m,n} ~,
\\
\bar{a}\ket{m,n}=a_{mn} \ket{m+1,n} ~,
\quad b\ket{m,n}=b_{mn} \ket{m,n+1} ~.
\end{array}
$$ 
By requiring that we have a $*$-representation we have also that
$$ a\ket{m,n}=a_{m-1,n}\ket{m-1,n}~, ~~ \bar{b}\ket{m,n}=b_{m,n-1}\ket{m,n-1} 
~,
$$ with the following recursion relations
$$ a_{m,n\pm 1} =q^{\pm 2} a_{m,n}~, \qquad b_{m \pm 1,n} = q^{\pm 2} b_{m,n}~,
\qquad b_{m,n} = q^2 a_{2n+1,m} ~.
$$ 
By solving them we get the explicit action of the generators on the basis of $\hil$:
\be\label{sigma}
\begin{array}{l}   t\ket{m,n}= q^{2m +4n+4} \ket{m,n},
\\    \bar a\ket{m,n}= (1-q^{2m+2})^{\frac{1}{2}}q^{m+2n+1}\ket{m+1,n},  \\
 a\ket{m,n}= (1-q^{2m})^{\frac{1}{2}}q^{m+2n}\ket{m-1,n}, 
\\  
b\ket{m,n}=(1-q^{4n+4})^{\frac{1}{2}}q^{2(m+n+2)}\ket{m,n+1},  \\ 
  \bar
b\ket{m,n}=(1-q^{4n})^{\frac{1}{2}}q^{2(m+n+1)}\ket{m,n-1}. 
\end{array}
\ee 
It is straightforward to check that all the defining relations
\eqref{s4} and
\eqref{sr4} are satisfied. Furthermore, 
The algebra generators are all trace class. Later one,  we shall need the trace of $t$
in this representation,
\be\label{classetraccia}
\Tr(t)=q^4 \sum_m q^{2m} \sum_n q^{4n}=\frac{q^4}{(1-q^2)(1-q^4)}. 
\ee
From the sequence of Schatten ideals in the algebra of compact 
operators one knows \cite{Sim79} that the norm closure of trace class operators gives the 
ideal of compact
operators $\ck$. Thus, the closure of $A(S_q^4)$ is the 
$C^*$-algebra
$\cc(S_q^4)=\ck \op \IC \II$. 

\subsection{The $K$-homology of $S_q^4$}

In the present situation we are
dealing with a $1$-summable Fredholm module $[\mu]\in K^0[\cc(S_q^4)]$. This is in 
contrast with the fact that the analogous element of $K_0(S^4)$ for the undeformed 
sphere is given by a
$4$-summable Fredholm module, being the fundamental class of  $S^4$.

For the Fredholm module $\mu := (\hil,\rho,\gamma)$ we have that the 
Hilbert space
is $\hil=\hil_\sigma \oplus \hil_\sigma$ with representation
$\rho=\sigma \oplus \beta$; here $\sigma$ is the representation of
$A(S^4_q)$ in \eqref{sigma} and $\beta$ given in \eqref{pi} is trivially extended to
$\hil_\sigma$. The grading operator is
\[
\gamma = \begin{pmatrix} 1&0
\\0&-1
\end{pmatrix} .
\] The corresponding Chern character $\chern^*(\mu)$ of the class of this Fredholm
module has  a component in degree $0$, $\chern^{0}(\mu) \in HC^0[A(S^{2n}_q)]$.
   From the general construction \cite{C94}, the element
$\chern^0(\mu_{\mathrm{ev}})$ is the trace
\[
\tau^1(x) := \Tr\left(\gamma \rho(x)\right) =
\Tr\left(\sigma(x) - \beta(x)\right).
\]
The operator $\sigma(x) - \beta(x)$ is always trace 
class. Obviously
$\tau^1(1)=0$. The higher degree parts of $\chern^*(\mu_{\mathrm{ev}})$ 
can again be
obtained via a periodicity operator.

\subsection{The charge and the rank} 
We are ready to compute the pairing giving the `topological charge'. Using \eqref{classetraccia} we find
\begin{align*}
\left\langle[\mu],[p]\right\rangle &:= \left\langle \chern^0(\mu),
\chern_0(p)\right\rangle = -q^{-4}(1-q^2)(1-q^4) ~\tau^1(t)  
\\  &= -q^{-4}(1-q^2)(1-q^4) \Tr(t)   \\
&= -q^{-4}(1-q^2)(1-q^4) q^{4}(1-q^2)^{-1}(1-q^4)^{-1} \\ &= - 1 ~.
\end{align*} 
This result also shows that the right $A(S^4_q)$-module
$p[A(S_q^4)^4]$ is not free. Indeed, any free module is represented in
$K_0[\cc(S_q^4)]$ by the idempotent $1$, and since
$\left\langle[\mu],[1]\right\rangle=0$, the evaluation of $[\mu]$ on any free
module always gives zero.

\bigskip 
We can extract the `trivial' element in the $K$-homology
$K^0[\cc(S_q^4)]$ of the quantum sphere $S_q^4$ and use it to measure 
the `rank' of
the idempotent
$p$. It corresponds to the trivial generator of the $K$-homology
$K_0(S^4)$ of the classical sphere
$S^4$. The latter (classical) generator is the image of the generator of the
$K$-homology of a point by the functorial map $K_*(\iota) : K_0(*) \to
K_0(S^{N})$, where $\iota : *
\hookrightarrow S^{N}$ is the inclusion of a point into the sphere. Now, the
quantum sphere
$S_q^4$ has just one `classical point', i.e. the $1$-dimensional representation
$\beta$ constructed in Sect.~\ref{se:pi}. The corresponding
$1$-summable Fredholm module $[\varepsilon]\in K^0[\cc(S_q^4)]$   is easily
described: the Hilbert space is $\IC$ with representation $\beta$; the grading
operator is
$\gamma=1$. Then the degree $0$ component $\chern^{0}(\varepsilon) \in
HC^0[A(S^{2n}_q)]$ of the corresponding Chern character is the trace 
given by the
representation itself:
\[
\tau^0(x) = \beta (x) ~,
\] 
and vanishes on all the generators whereas $\tau^0(1)=1$. \\ Not
surprisingly, the pairing with the class of the idempotent $p$ is
\[
\left\langle[\varepsilon],[p]\right\rangle := \tau^0(\chern_0(p)) =
\beta (2) = 2 ~.
\]

\appendix

\section{Connections on noncommutative vector bundles}\label{se:connections}

Let us suppose we have a differential calculus $(\omca=\op_p \oca{p}, \dd)$ for the algebra $\ca$. 
A  {\it connection} on the right finite projective $\ca$-module $\ce$ is a
$\IC$-linear map
\[
\nabla : \coca{p} \lra \coca{p+1} , 
\]
defined for any $p \geq 0$, and satisfying the Leibniz rule
\[
\nabla(\omega \rho) = (\nabla \omega) \rho + (-1)^{p} \omega \dd \rho  , 
~~\forall ~\omega \in \coca{p} , ~\rho \in \omca .  
\]
\noindent
A connection is completely determined by its restriction 
\[
\nabla : \ce  \to \coca{1}  ,
\]
 which satisfies
\[
\nabla (\eta a) = (\nabla \eta) a + \eta \ota \dd a  , 
~~\forall ~\eta \in \ce , ~a \in \ca ,  
\]
and which is extended by the Leibniz rule. It is again the latter property that implies the $\omca$-linearity of  
the composition, 
\[
\nabla^2 = \nabla \circ \nabla : \coca{p} \lra \coca{p+2} .
\]
The restriction of $\nabla^2$ to $\ce$ is the 
{\it curvature} 
\[
F : \ce  \to \coca{2} ,
\] 
of the connection. It is $\ca$-linear, $F(\eta a) =  F(\eta)a$ for any
$\eta\in\ce, a \in \ca$, and satisfies
\[
\nabla^2(\eta \ota \rho) =  F(\eta) \rho , ~~\forall ~\eta \in \ce , ~\rho \in \omca . 
\]
Thus, $ F\in \Hom_{\ca}(\ce, \coca{2} )$, the latter being the collection of (right) $\A$-linear homomorphisms of $\E$, with values in the two-forms $\Omega^2\A$.

In order to have the notion of a Bianchi identity we need some generalization. Let 
$\End_{\omca}(\comca)$ be the collection of all $\omca$-linear endomorphisms of $\comca$. It is an  algebra under composition. The curvature $ F$ can be thought of as an element of $\End_{\omca}(\comca)$. There is then a well-defined map
\begin{align} \label{def:conn-end}
[\nabla, ~\cdot~ ] & ~:~ \End_{\omca}(\comca) \lra \End_{\omca}(\comca) \nn \\
[\nabla, T] &:= \nabla \circ T - (-1)^{|T|}~ T \circ \nabla  .
\end{align}
where $|T|$ denotes the degree of $T$ with respect to the $\IZ^2$-grading of $\Omega\A$. 
It is easily checked that $[\nabla, ~\cdot~]$ is a graded derivation for the algebra $\End_{\omca}(\comca)$,  
\[
[\nabla,S \circ T]=[\nabla,S]\circ T + (-1)^{|S|} S \circ [\nabla,T].
\]
\begin{prop}\label{ubianchi}
The curvature $ F$ satisfies the  {\rm Bianchi identity},
\[
[\nabla,  F ] = 0 .  \label{ubia}
\]
\end{prop}
\begin{proof} 
Since $ F\in\End^0_{\omca}(\comca)$, the map $[\nabla,  F ]$ makes sense, and simply
$[\nabla,  F ] = \nabla \circ \nabla^2 - \nabla^2 \circ \nabla = \nabla^3 -
\nabla^3 = 0$. 
\end{proof}
\noindent

On the free module $\ce = \IC^N \otc \ca \simeq \ca^N$, a connection is given by the operator 
\[
\nabla_0 = \II \ot \dd : \IC^N \otc \oca{p} \lra \IC^N \otc \oca{p+1} .
\]
With the canonical identification $\IC^N \otc \omca = (\IC^N \otc \ca) \ota \omca ~\simeq~ (\omca)^N$, one thinks of $\nabla_0$ acting on $(\omca)^N$ as the
operator $\nabla_0 = (\dd, \dd, \cdots, \dd)$  ($N$-times).  
Take now any finite projective module with inclusion, $\lambda : \ce \to \ca^N$, identifying  $\ce$ as a direct summand of the free module $\ca^N$, and idempotent $p : \ca^N \to \ce$ which allows one to identify $\ce = p \ca^N$. Using these maps and their extensions to $\ce$-valued forms,  
a connection $\nabla_0$ on $\ce$ (called {\it Levi-Civita} or {\it Grassmann}) is the composition,
\[
\coca{p} ~\stackrel{\lambda}{\lra}~ \IC^N \otc \oca{p} 
~\stackrel{\II \ot  \dd}{\lra}~ \IC^N \otc \oca{p+1} 
~\stackrel{p}{\lra}~ \coca{p+1}, 
\]
that is  
\be\label{ugras}
\nabla_0 = p \circ (\II \ot  \dd ) \circ \lambda \simeq p \dd.    
\ee
All connections on $\ce$ constitute an affine space $C(\ce)$ modeled on the linear space $\Hom_{\ca}({\ce,\coca{1}})$.
Indeed, if $\nabla_1, \nabla_2$ are two connections on $\ce$, their difference is
$\ca$-linear,
\[
(\nabla_1 - \nabla_2)(\eta a) = ((\nabla_1 - \nabla_2)(\eta )) a , \quad \forall ~ \eta \in \ce , ~a \in \ca ,
\]
so that $\nabla_1 - \nabla_2 \in \Hom_{\ca}({\ce,\coca{1}})$. Thus, any connection can be written as 
\be
\nabla = p \dd + \alpha  , \label{uconn}
\ee
with $\alpha$ any element in $\Hom_{\ca}({\ce,\coca{1}})$.  The `matrix of $1$-forms'
$\alpha$ as in \eqref{uconn} is the {\it gauge potential}  of the connection $\nabla$ and the
corresponding curvature $ F$ is
\be
 F = p \dd p \dd p + p \dd \alpha + \alpha^2  . \label{ucurv} 
\ee

Next, let the algebra $\ca$ have an involution $^*$ extended to the whole of $\omca$ by the requirement 
$(\dd a)^* = \dd a^*$ for any $a\in\ca$.
A {\it Hermitian structure} on the module $\ce$ is a map $\hs{ \cdot}{\cdot} : \ce \ot \ce \to \ca$ with the properties
\begin{align} \label{hest}
&\hs{\eta a }{\xi } = a^* \hs{\xi}{\eta}  , \qquad
\hs{\eta}{\xi }^* = \hs{\xi}{\eta}  , \nn \\
&\hs{\eta}{\eta} \geq 0 ,  \qquad \hs{\eta}{\eta} = 0 \iff \eta = 0 , 
\end{align}
for any $\eta,\xi\in\ce$ and $a\in\ca$. We shall also require the Hermitian 
structure to be {\it self-dual}, that is every right $\ca$-module homomorphism 
$\phi: \E \to \ca$ is represented by an element of $\eta \in \E$, via the 
assignment $\phi(\cdot) = \hs{\eta}{\cdot}$, the latter having the correct 
properties by the first of \eqref{hest}. 
The Hermitian structure is naturally extended to a bilinear map from 
$\comca \times \comca$ to $\omca$ by
\be\label{uhext}
\hs{ \eta \ota \omega}{\xi \ota \rho }~ = (-1)^{|\eta| |\omega|} \omega^* \hs{\eta}{\xi}\rho , \quad \forall ~\eta, \xi \in \comca , ~\omega,\rho\in
\omca.
\ee

A connection $\nabla$ and a Hermitian structure $\hs{ \cdot}{\cdot}$ on $\ce$ are compatible if
\[
\hs{\nabla \eta}{\xi} + \hs{\eta}{\nabla \xi} = \dd \hs{\eta}{\xi}, \quad \forall ~\eta, \xi \in \ce .
\]
It follows directly from the Leibniz rule and \eqref{uhext} that this extends to
\[
\hs{\nabla \eta}{\xi} + (-1)^{|\eta|} \hs{\eta}{\nabla \xi} = \dd \hs{\eta}{\xi} , \quad \forall ~\eta, \xi \in \comca .
\]
We still use the symbol $C(\ce)$ to denote the space of compatible 
connections on $\ce$.

The collection $\End_{\ca}(\ce)$ of all $\ca$-linear endomorphisms of $\ce$, is an algebra with involution. Since we are taking a self-dual 
Hermitian structure, any $T \in 
\End_{\ca}(\ce)$ is adjointable, i.e. it admits an adjoint $T^*\in\End_{\ca}(\ce)$ such that
\[ \hs{T^* \eta}{\xi} = \hs{\eta}{T \xi} , \quad 
\forall ~\eta, \xi \in \ce . 
\]
The group $\cu(\ce)$ of unitary endomorphisms of $\ce$ is given by 
\[
\cu(\ce) := \{ u \in \End_{\ca}(\ce) ~|~ u u^* = u^* u = \id_\ce  \} .
\]
This group  plays the role of the {\it infinite dimensional group of gauge transformations}. 
It naturally acts on compatible connections by
\be\label{ugcon}
(u, \nabla) \mapsto \nabla^u := u^* \nabla u , \quad  \forall ~u\in \cu(\ce),  ~\nabla \in C(\ce) ,
\ee
where $u^*$ is really $u^* \ot \id_{\omca}$. Then the curvature transforms in a covariant way
\be
(u,  F) \mapsto  F^u = u^*  F u , \label{ugcur}
\ee
since, evidently, $ F^u = (\nabla^u)^2 = u^* \nabla u u^* \nabla u^* = u^* \nabla^2 u = u^*  F u$.

On the free 
module $\A^N$ there is a canonical Hermitian structure given by 
$\hs{\eta}{\xi} = \sum_{j=1}^N \eta_j^* \xi_j$, with 
$\eta = (\eta_1, \cdots, \eta_N)$ and $\eta = (\eta_1, \cdots, \eta_N)$ any two 
elements of $\A^N$. Under suitable regularity conditions on the algebra $\A$ all Hermitian structures on a given finite projective module $\E$ over $\A$ are isomorphic to each other and are obtained from the canonical structure on $\A^N$ by 
restriction \cite[II.1]{C94}. Moreover, if $\E=p \A^N$, then $p$ is self-adjoint,
$p=p^*$, with $p^*$ obtained by the composition of the involution $^*$ in the 
algebra $\A$ with the usual matrix transposition. 
The Grassmann connection \eqref{ugras} is easily seen to be compatible,  
\[ 
\dd \hs{\eta}{\xi} =  \hs{\nabla_0 \eta}{\xi} + \hs{\eta}{\nabla_0 \xi} .
\]
For a general connection \eqref{uconn}, compatibility reduces to
\[
\hs{ \alpha \eta}{\xi } + \hs{\eta}{\alpha \xi }= 0 , \quad  \forall ~\eta, \xi \in \ce , 
\]
which just says that the gauge potential is skew-hermitian,
\[
\alpha ^* = -\alpha . 
\]
Under the action \eqref{ugcon}, a gauge potential transforms in the
usual affine manner,
\be
(u, \alpha) \mapsto \alpha^u := u^* p \dd u + u^* \alpha u . \label{ugpot}
\ee
Indeed,  
$\nabla^u (\eta) = u^*( p \dd + \alpha) u \eta = u^* p \dd (u \eta) + u^* \alpha u \eta  
= u^* p u \dd \eta + u^* p (\dd u) \eta  + u^* \alpha u \eta  =  p \dd \eta + (u^* p \dd u + u^* \alpha u) \eta$
for any $\eta\in\ce$, yielding \eqref{ugpot} for the transformed potential.

\bigskip
Let $\End^s_{\omca}(\comca)$ denote the space of elements $T$ in $\End_{\omca}(\comca)$ which are skew-hermitian with respect to the Hermitian structure \eqref{uhext}, {i.e.}  satisfying
\be \label{def:skew}
\hs{T \eta}{\xi} + \hs{\eta}{T\xi}=0, \qquad \forall \; \eta,\xi\in \E.
\ee
A direct computation shows that the map $[\nabla,~\cdot~]$  
in \eqref{def:conn-end} restricts to the space $\End^s_{\omca}(\comca)$ as a derivation
\[
[\nabla, ~\cdot~ ]  ~:~ \End^s_{\omca}(\comca) \lra \End^s_{\omca}(\comca).
\]

The `tangent vectors' to the gauge group $\cU(\E)$ constitute the vector space of infinitesimal gauge transformations. For $X \in \End_{\Cinf(\S^4)}(\E)$ we define a family $\{u_t\}_{t \in \bR}$ of elements in $\cU(\E)$ by $u_t=1+t X +\cO(t^2)$, so that $X=(\partial u_t /\partial t)_{t=0}$. Unitarity of $u_t$ becomes $(1+t (X+ X^*)+ \cO(t^2))=1$. Taking the derivative with respect to $t$, at $t=0$, yields $X=-X^*$.
Thus, for $u_t$ to be a gauge transformation, $X$ should be a skew-hermitian endomorphisms of $\E$. In this way, we think of the real vector space $\End^s_\A(\E)$ as the collection of {\it infinitesimal gauge transformations}; its  complexification $\End^s_\A(\E)\ot_\bR \C$ can be identified with $\End_\A(\E)$. 

Infinitesimal gauge transformations act on connections in a natural way. Let the above family $u_t$ of gauge transformations  act on $\nabla$ as in \eqref{ugcon}. Since we have $(\partial (u_t \nabla u_t^* )/\partial t)_{t=0} = [\nabla, X]$, we conclude that an element $X \in \End^s_\A(\E)$ acts infinitesimally on a connection $\nabla$ by the addition of  $[\nabla,X]$,
\be\label{infgt}
(X, \nabla) \mapsto \nabla^X = \nabla + t [\nabla,X] + \cO(t^2) , \quad \forall ~X \in \End^s_\A(\E),  ~\nabla \in C(\ce). 
\ee
As a consequence, for the  transformed curvature one finds 
\[
(X,  F) \mapsto F^X = F + t [F,X] + \cO(t^2), 
\]
since $F^X = (\nabla + t [\nabla,X])\circ (\nabla + t [\nabla,X]) = \nabla^2 + t [\nabla^2, X] + \cO(t^2)$.

\section{The Chern character of projections}\label{se:cc}

Let $C_*(\A)$ be the complex consisting of cycles over the algebra $\A$, that is $C_n(\A):=\A^{\otimes(n+1)}$ in degree $n$. On this complex there is defined the Hochschild boundary operator $b:C_n(\A)\to C_{n-1}(\A)$ given
by
\[
b(a_0 \ot a_1 \ot \cdots \ot a_n) :=
\sum_{j=0}^{n-1} (-1)^j a_0 \ot \cdots \ot a_j a_{j+1} \ot \cdots \ot a_n 
+ (-1)^{n} a_n a_0 \ot a_1 \ot \cdots \ot a_{n-1} \, .
\]
It is easy to prove that $b^2=0$.
The {\it Hochschild homology} $HH_{*}(\ca)$ of the algebra $\ca$ is
the homology of this complex. There is a second operator which increases the degree
$ B:\cha{n}\to \cha{n+1}$ written as 
$ B =  B_0 N$, with 
\begin{align*}
B_0 (a_0 \ot a_1 \ot \cdots \ot a_n ) &:= \II \ot a_0 \ot a_1 \ot \cdots \ot a_n
\\
N (a_0 \ot a_1 \ot \cdots \ot a_n ) &:= \frac{1}{n + 1} \sum_{j=0}^n 
(-1)^{nj} a_j \ot a_{j+1} \ot \cdots a_n \ot a_0 \ot \cdots \ot a_{j-1} \,  . 
\end{align*}
with the obvious cyclic identification $n+1 = 0$.
Again it is straightforward to check that $B^2=0$ and that $b B + B b = 0$;  thus $(b+B)^2=0$. By putting together these two operators, one gets a bi-complex
$(C_*(\ca), b, B)$ with $\cha{p-q}$ in bi-degree $p,q$. 
The {\it cyclic homology} $HC_{*}(\ca)$ of the algebra $\ca$ is
the homology of the total complex $(CC(\ca), b+B)$,  
whose $n$-th term is given by 
\[
\ccha{n}  := \op_{p+q=n} \cha{p-q} = \op_{0\leq q
\leq  [n/2]} \cha{2n-q}.
\]

\bigskip
For $e \in K_0(\A)$, the Chern character $\chern_*(e) = \sum_{k\geq 0} \chern_k(e)$ is the even cyclic cycle in $CC_*(\A)$, $(b+B)\chern_*(e)=0$, defined by the following formul{\ae}. For $k=0$,
\[
\chern_0(e):= \tr(e),
\]
whereas for $k\neq 0$
\[
\chern_k(e):=(-1)^k \frac{(2k)!}{k!}  \sum (e_{i_0 i_1}-\frac{1}{2}\delta_{i_0 i_1}) \otimes {e_{i_1 i_2} \otimes 
e_{i_2 i_3} \otimes \cdots \otimes e_{i_{2k} i_0} }.
\]

If  $(\A,\cH,D,\gamma)$ is a spectral triple, the components of the Chern character are represented as bounded operators on  the Hilbert space $\ch$ by explicit formul{\ae},
\[
\pi_D(\chern_k(e))= (-1)^k \frac{(2k)!}{k!}  \sum (\pi(e_{i_0 i_1})-\frac{1}{2}\delta_{i_0 i_1}) 
[D, {\pi(e_{i_1 i_2})] 
\cdots [D,\pi(e_{i_{2k} i_0})] }.
\]
These operators are used in noncommutative index theorems, a very simple instance of this being given by \eqref{topbasic} when computing a topological number.

\bibliographystyle{amsalpha}

\end{document}